\documentclass[12pt]{article}

\usepackage{fix-cm}
\usepackage{mathtools}
\usepackage{enumerate}
\usepackage{amssymb,mathrsfs,amsthm,multirow}
\usepackage{booktabs}
\usepackage{siunitx}
\usepackage{graphicx}
\usepackage{subfig}
\usepackage{xcolor}
\usepackage{graphicx}
\usepackage{chngcntr}
\usepackage{amsmath,cases,multirow,url}
\usepackage{booktabs,bm}

\usepackage{geometry}
\usepackage{amssymb,mathrsfs,amsthm,multirow}

\usepackage{mathtools}
\usepackage{indentfirst}    
\usepackage{xcolor}
\usepackage{cite}
\usepackage{graphicx}
\usepackage{chngcntr}
\usepackage[inline]{enumitem}
\usepackage[ruled,linesnumbered]{algorithm2e}

\usepackage{tikz}
\usetikzlibrary{positioning,shapes}
\usetikzlibrary{matrix}
\usetikzlibrary{graphs}
\usetikzlibrary{chains}

\newcommand{\xiaosi}{\fontsize{12pt}{12pt}\selectfont}      

\newcommand{\setd}{{ d \kern -.15em l}}
\newcommand{\hatsetd}{ d \hat{\kern -.15em l }}
\newcommand{\dd}{\mathsf {d\kern -0.07em l}} 

\newcommand{\bgeqn}{\begin{eqnarray}}
\newcommand{\edeqn}{\end{eqnarray}}
\newcommand{\bgeq}{\begin{eqnarray*}}
\newcommand{\edeq}{\end{eqnarray*}}
\newcommand{\bec}{\begin{center}}
\newcommand{\enc}{\end{center}}

\newtheorem{lemma}{Lemma}

\newtheorem{assumption}{Assumption}
\usepackage{color,colortbl,xcolor}

\begin{document}
\bibliographystyle{unsrt}

\title{
A Robo-Advisor System: expected utility modeling via pairwise comparisons
\footnote{This research was supported by National Key R\&D Program of China under No. 2022YFA1004000 and National Natural Science Foundation of China under Grant Number 11991023 and 12371324.}}

\author{Bo Chen, Jia Liu\footnote{Emails:chenbo1211@stu.xjtu.edu.cn , jialiu@xjtu.edu.cn} \\
\footnotesize
School of Mathematics and Statistics, Xi'an Jiaotong University, 710049, Xi'an, P. R. China
}
\date{}

\maketitle 

\begin{abstract}\xiaosi
We introduce a robo-advisor system that recommends customized investment portfolios to users using an expected utility model elicited from pairwise comparison questionnaires. 
The robo-advisor system comprises three fundamental components. First, we employ a static preference questionnaire approach to generate questionnaires consisting of pairwise item comparisons.
Next, we design three optimization-based preference elicitation approaches to estimate the nominal utility function pessimistically, optimistically, and 
{\color{black}neutrally.}
Finally, we compute portfolios based on the nominal utility using an expected utility maximization optimization model.
We conduct a series of numerical tests on a simulated user and 
{\color{black}a number of} human users to evaluate the efficiency of the proposed model.
\end{abstract}


\section{Introduction}
\setcounter{page}{1}
\pagenumbering{arabic}

Robo-advisors are a class of financial advisers providing financial advice and investment management online with moderate to minimal human intervention, 
which provide 
personalized financial advice based on mathematical rules or algorithms. 
A robo-advisor system provides personalized investment advice to 
users with different ages, investment goals, investment periods, risk preferences, and even preferred investment methods. 
As a result, the investment recommendations from these services vary from one person to another. 
Due to their ability to offer 
inclusive and personalized financial services to a broader group of users
at a low cost, robo-advisors have rapidly evolved and found applications in various areas, such as brokerage \cite{Acunto2019} or personal advisor services  \cite{Rossi2020}. 
In some 
emerging frontiers, robo-advisors are also being applied to manage debt repayment \cite{Chak2022} and to assist in determining optimal consumption and spending by, for instance, 
utilizing big data to construct benchmarks based on the behavior of peers \cite{Acunto2023}.

Broadly speaking, the portfolio recommendation methodology in research of robo-advisors can be divided into two subjects: the model-free approach
and the model-based approach.  
The model-free approaches are commonly based on machine learning methods, such as 
unsupervised learning or reinforcement learning.
For instance, some robo-advisors recommend investment portfolios or investment items by matching the characteristics 
of the user with
historical users.
One effective technique is collaborative filtering, which recommends financial products or investment strategies that are most similar 
to the historical choices of the user, or recommendations 
based on choices from similar historical users of the robo-advisor \cite{Sarwar2001, Zhao2010}.  
Reinforcement learning approaches gradually improve the fitness of the recommendations by continuously interacting with the user and collecting feedback~\cite{Bourdache2019, Alsabah2021, wjx}.


This paper focuses on the model-based approaches, which usually set a modern financial optimization model first, then design an interaction mechanism with the user consistent with the base model. 
On one side, the robo-advisor learns the user's preferences, represented by some parameters in the optimization model, through some elicitation models or algorithms during the interaction \cite{toubia2003fast}. 
On the other side, the robo-advisor provides the most suitable portfolio to the user based on the model with {\color{black}the} elicited parameters 
\cite{bertsimas2013learning,Cui2022}. 
We can apply different modern financial optimization models as the guidance model of the user. 
Alsabah et al.~\cite{Alsabah2021} propose a reinforcement learning framework for robo-advising, 
in which they consider several mean-risk 
models characterized by risk-aversion parameters as the user's true criteria for portfolio selection.
{\color{black}
In a similar vein, Dong et al.~\cite{Dong2022} 
employ a bi-level mean-variance framework 
in their proposed robo-advisor 
and provide a closed-form 
for the risk-aversion parameter.} 
However, there exist non-monotonicity and time-inconsistency issues
in the dynamic mean-variance optimization, and Cui et al.~\cite{Cui2022} propose a novel dynamic asset allocation framework based on {\color{black}a family of extended mean-risk models in which the target wealth and risk-aversion parameter of the user 
can be elicited by the robo-advisor system}.
Inverse optimization methods are also applied to 
elicit 
the risk preference from the user's historical
portfolios under a mean-variance portfolio allocation framework~\cite{Yu2023}. 


Different from the mean-risk framework which seeks an optimal trade-off between risk and return, the famous expected utility maximization framework proposed 
by von Neumann and Morgenstern~\cite{Von1947} searches for the investment portfolio that maximizes the utility function of a rational user. 
They 
state that there exists a utility function 
$u:\mathbb{R}\rightarrow\mathbb{R}$, so that
any set of preferences a user may hold between risky outcomes (represented by random variables $W$ and $Y$) can be formulated by an expected utility measure if the preferences 
adhere to certain reasonable axioms, i.e., $\mathbb{E}[u(W)]\geq\mathbb{E}[u(Y)]$ if and only if the user prefers $W$ to $Y$. 
However, despite its mathematical rigor and flexibility, expected utility maximization has not been studied by the robo-advisor community due to its non-parametric nature and the difficulties in its preference elicitation. 


Consequently, the major difficulty of applying the expected utility framework in a robo-advisor system is the estimation/learning of the user's utility function. 
According to the information observed by the system, which characterizes users' preferences, there are three major approaches for estimating utility functions in the preference learning area. 
The first one employs {\color{black}pairwise comparisons}, where the system provides some items/lotteries in pairs to the user and collects the user’s choices between the items \cite{eggers2021choice}. 
Then, a conjoint analysis method or preference robust optimization method can be applied to identify the value of the utility function at a discrete set of points for elicitation~\cite{toubia2003fast,Delage2015}. 
{\color{black}Alternatively, 
the preference graph from pairwise comparisons is also an extensively discussed approach.
For instance, \v{C}aklovi\'c and Kurdija \cite{vcaklovic2017universal} show its applications to Eurovision Song Contest and Cookie type election, while Csat\'o and T\'oth \cite{csato2020university} present a methodology to rank universities on the basis of the lists of programs the students applied for.}
The second approach collects the user's ratings on certain products and utilizes these ratings to construct a utility function~\cite{MS2021}. 
The last approach is to observe the decisions of the user in historical decision-making problems and elicit the utility by {\color{black} inverse optimization
model~\cite{baak2023preference}} 
or inverse reinforcement learning~\cite{inverse-RL2023}. 
The most widely used approach is the pairwise comparison (choice-based) method, which is easy to implement 
in the interaction 
with users, as surveyed in 
\cite{eggers2021choice,Clemen2013}. 
For instance, it can be observed in many risk tolerance assessment surveys used by financial advisers, as these typically involve questions such as the following, from Grable and Lytton~\cite{Grable1999}: 
\begin{quote}
    You are on a TV game show and can choose one of the following. Which would you take?

    (a)\quad \$1,000 in cash,

    (b)\quad A 50\% chance at winning \$5,000, 

    (c)\quad A 25\% chance at winning \$10,000, 

    (d)\quad A 5\% chance at winning \$100,000. 
\end{quote}
From the answers/choices of the user to the items, we can estimate the parameters in the utility function {\color{black}with} the pre-set form \cite{Jung2018a, Jung2018b}.
{\color{black}Most applications of expected utility in portfolio selection use parameterized expected utility functions.} 
Correspondingly, traditional elicitation methods with pairwise comparison data focus on particular types of utility functions. For example, {\color{black}the conjoint analysis~\cite{Rao2014} applies a 
linear utility function to measure the user' preferences for different attributes of the products and the robo-advisor system only needs to learn the loading factors of the linear utility function \cite{wjx}.
However, the structure of the utility function of a user
cannot be pre-specified in more broad application cases, which might potentially be linear, quadratic, or exponential}.
Designing an optimization-based method to elicit a non-parametric utility function (without pre-specifying the parameterized utility function form) consistent with the choices of pairwise comparisons is still a challenge in this research area. 
It is worth noting that there is another stream of research to handle the ambiguity of utility in the elicitation process of preference by using the preference robust approach \cite{Delage2015,hu2017optimization,guo2021statistical,guo2023utility}. This approach makes robust decisions based on the worst-case utility function from the ambiguity set constructed by pairwise comparisons or moments information.

Another important issue in the pairwise comparison approach is the generation method of the corresponding questionnaires, which directly impacts the efficiency 
of the preference elicitation. 
Particularly in the robo-advisor applications, 
applying some mathematical models or algorithms to generate questionnaires is more suitable for 
minimizing human intervention.
Christakopoulou et al.~\cite{Christakopoulou2016} 
introduce a bandit-based method for creating a dynamic questionnaire with relative questions in an online setting, in which an algorithm is designed to pick two items for relative questions. 
{\color{black}To select $K$ pairwise items all at once and 
generate a static preference questionnaire consisting of $K$ relative questions, Sepliarskaia et al.~\cite{SPQ2018} propose an 
algorithm based on the latent factor model.
Moreover, Sz\'adoczki et al. \cite{szadoczki2022filling} study the decision problems where the set of pairwise comparisons can be chosen completely before the decision making process without any further prior information.}
Armbruster and Delage~\cite{Delage2015} propose the random utility split and random relative utility split 
methods to determine the pre-set 
parameters in simple questions
in the form
“Do you prefer a certain return of $r_2$ or a lottery where the return will be $r_3$ with probability $p$ and $r_1$ with probability $1-p$?”
However, preference elicitation and questionnaire generation are separately designed and studied in traditional research. 
To the best of our knowledge, there is not a unified robo-advisor model in the expected utility framework. In this paper, we will build a basic framework for non-parametric robo-advisor systems and study the model design and solution methods in different components.


\begin{figure}[htbp]
	\centering 
	\includegraphics[scale=0.5]{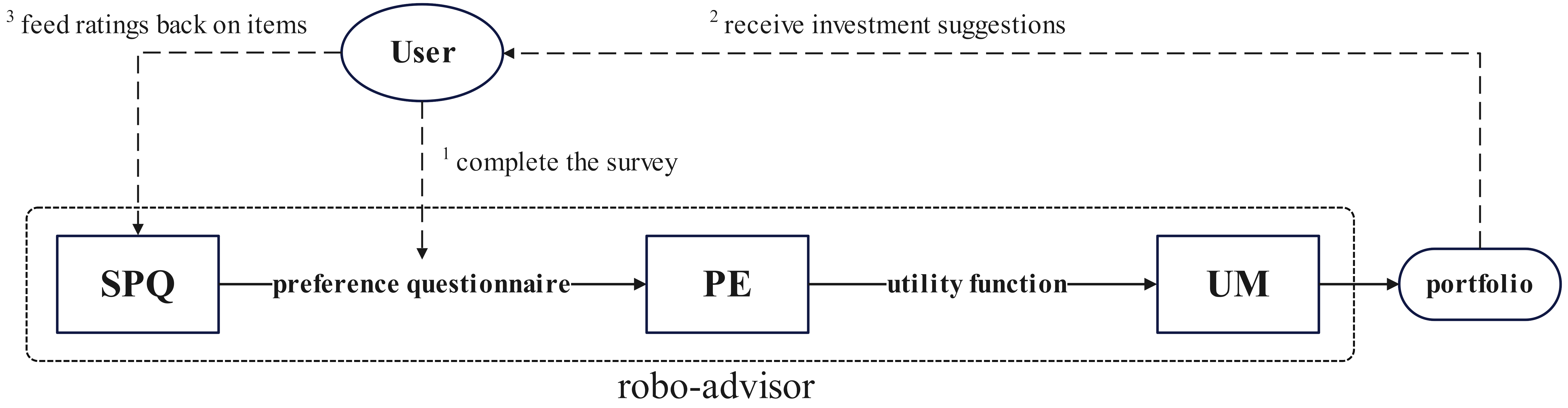}
        \caption{Flowchart of the robo-advisor system}
        \label{robo-advisor flowchart}
\end{figure}

In this paper, we design a new unified framework of a robo-advisor system 
constituted by three key components: 
static preference questionnaire (SPQ), preference elicitation (PE), utility maximization (UM), as shown in Figure~\ref{robo-advisor flowchart}. 
The SPQ component generates preference questionnaires comprising $K$ pairs of 
items for the user to make binary choices, thus forming corresponding pairwise comparisons that are utilized in the subsequent preference elicitation procedure.
Specifically, we formulate the task of selecting $K$ pairs of items from a preset item set as an optimization problem and employ a greedy iteration algorithm to solve it. 
Then, in the preference elicitation process, we adopt a non-parametric expected utility framework and develop the user's utility function by 
piecewise linear approximation for portfolio recommendation.
In particular, we construct an ambiguity set consisting of all utility functions consistent with the user's pairwise comparisons and determine three different nominal utility functions from the set under pessimistic, optimistic, and neutral estimations, using optimization methods.
The final recommended portfolio is obtained by maximizing the user's utility on historical return data, with the nominal utility functions derived from preference elicitation. We summarize the main notations used in the paper in Table \ref{tab:notations}.

\begin{table}[htbp]
\centering
\caption{Notations}\label{tab:notations}
\resizebox{\textwidth}{!}{
\begin{tabular}{l l }
\toprule
\multicolumn{2}{l}{\textbf{Input parameters for SPQ:}} \\
$\mathbb{I}^n=\{I_1,\ldots, I_n\}$ &  set of $n$ items\\
$I_i$ & the $i$-th item in $\mathbb{I}^n$\\
$R$ & historical rating matrix provided by historical users\\
$r_{ui}$ & rating of item $I_i$ rated by historical user $u$ \\
$K$ & number of questions in SPQ \\
{\textbf{Decision variables in SPQ:}} \\
$p_i$  &  latent factor vector of item $I_i$  \\
$b_i$ & bias of item $I_i$\\
$q_u$ & latent factor vector of user $u$ \\
$s_u$ & bias of user $u$ \\
$\hat{a}_{uij}$ & true relative preference of item $I_i$ over item $I_j$ of user $u$\\
$\tilde{a}_{uij}$ & predicted relative preference of item $I_i$ over item $I_j$ of user $u$\\
$\mathbb{B}$ & set of selected item pairs in SPQ\\
$(W_k, Y_k)$ & the $k$-th pair of items in $\mathbb{B}$\\
\multicolumn{2}{l}{\textbf{Input parameters for preference elicitation:}} \\
$\overline{\mathbb{Y}}$& ordered set of breakpoints\\
$\bar{y}_j$ & the $j$-th smallest entry in $\overline{\mathbb{Y}}$\\
$\eta=\left[I_1, I_2, \ldots, I_n\right]^{\top}$ & random vector constituted by all items in $\mathbb{I}^n$
\\
$\eta_i$ & the realization of $\eta$ under the $i$-th scenario 
\\
$p_{i}$ & probability of $\eta=\eta_i$\\
$N$ & cardinality of the set $\overline{\mathbb{Y}}$\\
$M$ & number of scenarios of $\eta$ \\
\multicolumn{2}{l}{\textbf{Decision variables in preference elicitation:}} \\
$u_N$ & piecewise linear utility function with $N$ pieces 
\\
$\alpha_j$ & utility value at the $j$-th breakpoint in $\overline{\mathbb{Y}}$ \\
$\beta_j$ & slope of the $j$-th linear segment 
in $u_N$\\
\multicolumn{2}{l}{\textbf{Input parameters for portfolio selection:}} \\
$u^X, X=P,O,N$ & nominal utility function under pessimistic(P), optimistic(O),\\
& or neutral(N) estimation\\
$\xi$ & random return vector of $S$ risky assets\\
$\xi_t$  & realization of return rate 
on the $t$-th historical day\\
$W_0$ & initial investment budget amount\\
\multicolumn{2}{l}{\textbf{Decision variables in portfolio selection:}} \\
$x=\left[x_0,x_1,x_2,\ldots,x_S\right]$ & {\color{black}portfolio on a risk-free asset and $S$ risky assets}\\



\bottomrule
\end{tabular}}
\end{table}

The major contributions of this paper can be summarized in the following: 
\begin{itemize}
[itemsep=2pt]
    \item We propose the first robo-advisor system in 
    the non-parametric expected utility framework. 
    \item We propose new optimization-based preference elicitation methods from pairwise comparison data.
    \item We employ the static preference questionnaire to generate pairwise comparison questions 
    which can achieve effective preference elicitation with fewer questions. 
    \item We conduct numerical tests to validate the robo-advisor system and particularly verify it with 
    human participants.  
    
\end{itemize}

This paper is organized 
as follows. 
In Section \ref{sec-spq}, we introduce the SPQ method to generate questionnaires with pairwise comparison questions 
provided to users. 
We then propose three preference elicitation models from the user's choice data on the pairwise comparisons in Section \ref{sec-pre-elic}. 
In Section \ref{sec-portfolio-opt}, we apply a utility maximization model 
to compute the recommended investment portfolios.
Section \ref{sec-numerical} includes various numerical tests on simulated users and 
human users to demonstrate the performance of the robo-advisor. Section \ref{sec-conclusion} concludes.

\section{Static preference questionnaire}\label{sec-spq}

When a new user visits the robo-advisor, the core of the user interaction is to provide the user with some pairwise comparison questions between pre-set items. 
An item is set as a lottery including several realizations of returns with corresponding probabilities, and each question in the questionnaire provides a pair of items $(I_i, I_j)$ to the user for pairwise comparison, 
{\color{black}such as: }
\begin{quote}
    Suppose you get a chance to win a prize. Choose the preferred one between the following.

    (a)\quad A 1\% chance at winning \text{\textyen}500,000,

    (b)\quad A 50\% chance at winning \text{\textyen}5,000.
\end{quote}

We first pre-set an item set and use the static preference questionnaire algorithm to pick certain numbers of pairs of items from the item set. 
Denote the item set with $n$ items as $\mathbb{I}^n:=\{I_1,I_2,\ldots,I_n\}$, where $I_i$ is a non-negative valued discrete random variable corresponding to the $i$-th item. 
We require the items in $\mathbb{I}^n$ to be distinct. 
An example set with 10 items $\mathbb{I}^{10}$ is shown in Table~\ref{10-item set}. 

\begin{table}[htbp]
\renewcommand{\arraystretch}{0.9} 
\setlength{\extrarowheight}{3.5pt}
  \centering
  \caption{An instance of item set $\mathbb{I}^{10}$}\label{10-item set}
    \resizebox{\textwidth}{!}{
    \begin{tabular}{l l l}
    \toprule
    Item 
    & Description & 
    Probability representation \\
    \midrule
    1 & \text{\textyen}800 in cash & $\mathbb{P}(I_1=800)=1$ \\
    2 & An 80\% chance at winning \text{\textyen}1,000 
    & $\mathbb{P}(I_2=1000)=0.8,\ \mathbb{P}(I_2=0)=0.2$ \\
    3 & A 50\% chance at winning \text{\textyen}5,000 
    & $\mathbb{P}(I_3=5000)=0.5,\ \mathbb{P}(I_3=0)=0.5$ \\
    4 & A 25\% chance at winning \text{\textyen}10,000 
    & $\mathbb{P}(I_4=10000)=0.25,\ \mathbb{P}(I_4=0)=0.75$ \\
    5 & A 5\% chance at winning \text{\textyen}100,000 
    & $\mathbb{P}(I_5=100000)=0.05,\ \mathbb{P}(I_5=0)=0.95$ \\
    6 & A 1\% chance at winning \text{\textyen}500,000 
    & $\mathbb{P}(I_6=500000)=0.01,\ \mathbb{P}(I_6=0)=0.99$ \\
    7 & A 0.1\% chance at winning \text{\textyen}1,000,000 
    & $\mathbb{P}(I_7=1000000)=0.001,\ \mathbb{P}(I_7=0)=0.999$ \\
    8 & \parbox[t]{8cm}{A 50\% chance at winning \text{\textyen}1,000 \\ [-0.1ex] and a 10\% chance at winning \text{\textyen}10,000}
    & \parbox[t]{8.2cm}{$\mathbb{P}(I_8=1000)=0.5,\ \mathbb{P}(I_8=10000)=0.1,\\ [-0.1ex] \mathbb{P}(I_8=0)=0.4$}\\
    9 & \parbox[t]{8cm}{A 20\% chance at winning \text{\textyen}10,000\\[-0.1ex]
    and a 1\% chance at winning \text{\textyen}100,000} 
    & \parbox[t]{8.2cm}{$\mathbb{P}(I_9=10000)=0.2,\ \mathbb{P}(I_9=100000)=0.01,\\ [-0.1ex] \mathbb{P}(I_9=0)=0.79$} \\
    10 & \parbox[t]{8cm}{A 50\% chance at winning \text{\textyen}2,000 \\[-0.1ex]
    and a 2\% chance at winning \text{\textyen}200,000} 
    & \parbox[t]{8.3cm}{$\mathbb{P}(I_{10}=2000)=0.5,\ \mathbb{P}(I_{10}=200000)=0.02,\\[-0.1ex] \mathbb{P}(I_{10}=0)=0.48$} \\
    \bottomrule
    \end{tabular}}
\end{table}

\begin{figure}[htbp]
    \centering
    \includegraphics[scale=0.45]{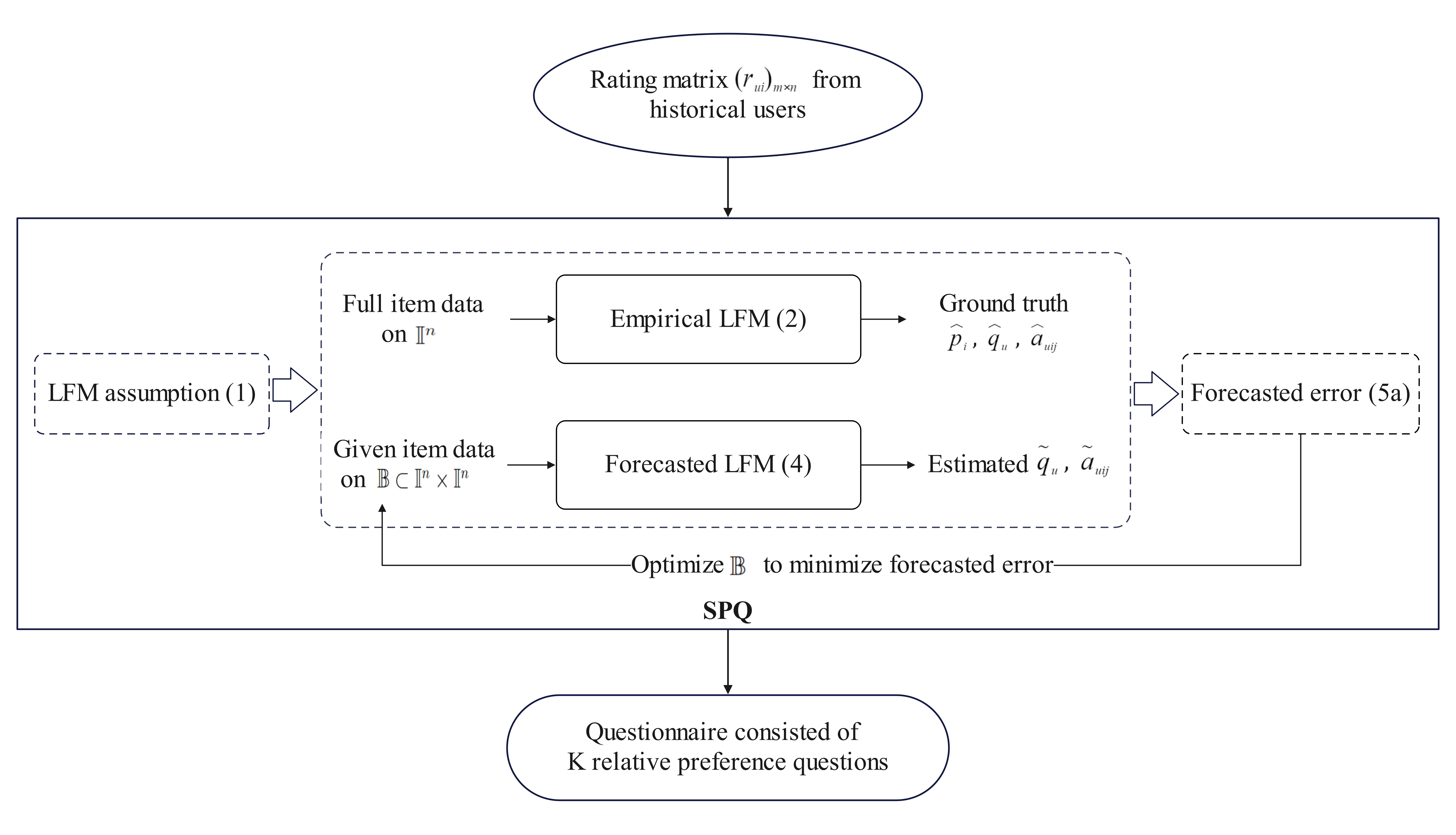}
    \caption{Flowchart of the static preference questionnaire generation method}
    \label{fig:SPQ}
\end{figure}

We then apply the static preference questionnaire (SPQ) generation method, motivated by \cite{Anava2015}, to select some pairs of items for the questionnaires. 
The basic idea of the SPQ method is to collect some empirical rating scores from historical users on a full set of items, and then select several pairs of items that 
reflect the users' preference to the 
{\color{black}highest}
extent. 
The SPQ method assumes a latent factor model and fits the empirical latent factors with the full set of data. 
Then we apply an optimization model to find some pairs of items in which the forecasted values of rating scores are mostly close to the empirical values. 
A brief flowchart of the SPQ method is shown in Figure~\ref{fig:SPQ}.

\subsection{Assumptions of SPQ}\label{sec-spq-aasump}
{\color{black}
Suppose there are $m$ historical users in the robo-advisor system, from whom we have collected rating scores for items in the pre-set item set. Several feasible approaches exist to obtain these scores.
One straightforward method is providing all items from the pre-set item set, as shown in Table~\ref{10-item set}, to historical users. Each historical user is then asked to assign a score to each item, usually ranging from 0 to 10. However, a potential issue arises with missing data, wherein some historical users may overlook rating certain items. To address this issue, we can complete or infer the missing data using known scores. Alternatively, in subsequent analyses, we can just eliminate the regression error terms corresponding to the missing data in \eqref{LFM_solve}.
Another approach is to utilize existing models to estimate historical user ratings for items based on their historical pairwise comparison choices, such as conjoint analysis~\cite{toubia2003fast}, preference graphs~\cite{vcaklovic2017universal,csato2020university},
random utility theory~\cite{yellott1977relationship, liu2024robust} or
Sonneborn–Berger score in the Swiss-system~\cite{sziklai2022efficacy}.
}

We denote the historical rating matrix as $R\in \mathbb{R}^{m\times n}$, where each entry $r_{ui}$ denotes the rating of item $I_i$ given by the historical user $u$, $i=1, \ldots, n$, $u=1, \ldots, m$.
For each item in $\mathbb{I}^n$, we aim to identify a $v$-dimensional latent factor vector to represent the characteristics of the item, which is obtained only from historical rating data. 
We assume a Latent Factor Model (LFM) \cite{Mnih2007} to explain the rating: 
\begin{equation}\label{LFM_model}
    r_{ui} = \mu + b_i + s_u + p_i^{\top} q_u + \varepsilon_{ui},
\end{equation}
where $\mu$ presents the global mean of all ratings in $R$, 
{\color{black}$p_i\in \mathbb{R}^v$ is the latent factor vector of item $I_i$ and $b_i\in\mathbb{R}$ is the bias of item $I_i$, $i=1, \ldots, n$, $q_u\in\mathbb{R}^v$ is the latent factor of user $u$ and $s_u\in\mathbb{R}$ is the bias of user $u$, $u=1, \ldots, m$, }
$\varepsilon_{ui}$ is the random rating noise. 
We then employ the 
Least Squares 
{\color{black}(LS)}
method to 
{\color{black}estimate}
the LFM model by solving the following optimization problem: 
\begin{equation}\label{LFM_solve}
    \min\limits_{p, q, b, s} \sum_{\substack{i=1,\ldots,n \\ u=1,\ldots,m}} (r_{ui}-p_i^{\top}q_u-b_i-s_u-\mu)^2 + \lambda_1\sum_{u=1}^{m} (\left\| q_u\right\|_2^2+\left\| s_u\right\|_2^2) + \lambda_2\sum_{i=1}^{n} (\left\| p_i\right\|_2^2+\left\| b_i\right\|_2^2). 
\end{equation}
{\color{black}In \eqref{LFM_solve}, we add a $L_2$-regularization term for each variable, totaling $2(m+n)$ terms. These terms constitute a penalty for hypothesis complexity, contributing to the prevention of overfitting \cite{shafieezadeh2019regularization}.}
We denote the optimal solutions of 
problem~\eqref{LFM_solve} as $\hat{p}_i,\hat{b}_i, i=1,\ldots, n$, $\hat{q}_u, \hat{s}_u, u=1,\ldots,m$. 
The larger the value of the product ${\hat{p}_i}^{\top}\hat{q}_u$, the more user $u$ prefers item $I_i$. 

\subsection{SPQ}
We then aim at finding a set $\mathbb{B}=\{(I_k^1, I_k^2), k=1,\ldots, K \mid I_k^1, I_k^2 \in \mathbb{I}^n \}\subset \mathbb{I}^n\times \mathbb{I}^n$, which represents the questionnaire constituted with $K$ pairwise comparison questions that we provide to the user. 

The selection principle of the SPQ method is to minimize the error between the empirical LFM model trained by all items' data and the forecasted LFM model estimated by the data of selected items. 
As we aim at finding pairwise items, we consider the relative preference between any two items, which is defined as the difference in the scores of the two items. 
Different from directly comparing the scores in \cite{Anava2015}, the relative preference we introduce can eliminate the effect of items to the greatest extent. 
Compared to \cite{SPQ2018}, which should collect extra binary choice data, our approach only relies on historical scoring data. 
By \eqref{LFM_model}, the relative preference on $(I_i, I_j)$ estimated by all items' data is 
\begin{equation}\label{true_preference}
    \hat{a}_{uij} = (\hat{p}_i-\hat{p}_j)^{\top}\hat{q}_u + (\hat{b}_i-\hat{b}_j) + \varepsilon_{ui}-\varepsilon_{uj}, \ u=1,\ldots,m,\ i,j=1,\ldots,n, 
\end{equation}
where $\hat{q}_u$, $u=1,\ldots,m$, $\hat{p}_i$, $\hat{b}_i$, $i=1,\ldots,n$ are the 
estimated values by solving \eqref{LFM_solve}. 

As for the forecasted LFM model, we estimate the forecasted relative preference on $(I_i, I_j)$ of user $u$ with the data of selected items in $\mathbb{B}$, by solving the following least 
{\color{black}squares}
optimization problem: 
\begin{equation}\label{get_qu}
    \min\limits_{q_u} \sum_{(I_i, I_j)\in \mathbb{B}}\left[a^*_{uij}- (\hat{p}_i-\hat{p}_j)^{\top}q_u-(\hat{b}_i-\hat{b}_j) \right]^2,\ u=1,\ldots,m, 
\end{equation}
where $\hat{p}_i-\hat{p}_j$ represents the relative difference between $I_i$ and $I_j$ and $\hat{b}_i-\hat{b}_j$ is the difference of bias. 
{\color{black}$a^*_{uij}:=r_{ui}-r_{uj}$ is the ground truth of the relative preference of user $u$, i.e., the difference between original score data $r_{ui}$ and $r_{uj}$.}
Noticing that according to LFM model assumed in \eqref{LFM_model}, $a^*_{uij}, u=1,\ldots,m$ also holds the form as \eqref{true_preference}. 
We denote the optimal solution of \eqref{get_qu} as $\tilde{q}_u, u=1,\ldots,m$, then we can use $\tilde{q}_u$ to forecast the relative preference on all other pairwise items by
\begin{equation*}
    \tilde{a}_{uij}=(\hat{p}_i-\hat{p}_j)^{\top}\tilde{q}_u+(\hat{b}_i-\hat{b}_j)+\varepsilon_{ui}-\varepsilon_{uj}, \ u=1,\ldots,m,\ i,j=1,\ldots,n. 
\end{equation*}

Then we minimize the 
{\color{black}squared}
error between $\hat{a}_{uij}$ and $\tilde{a}_{uij}$ on all pairwise items to optimize the selection of pairwise items in $\mathbb{B}$, which can be formulated as the following optimization model:
\begin{subequations}\label{SPQ}
\begin{align}
\min\limits_{\mathbb{B}}\ & \mathbb{E}\left[ \sum_{(I_i, I_j)\in \mathbb{I}^n\times \mathbb{I}^n }\left( \tilde{a}_{uij}-\hat{a}_{uij}\right)^2\right]\\
\text{s.t.}\ \ & \mathbb{B}\subset \mathbb{I}^n\times\mathbb{I}^n,\ \left|\mathbb{B}\right|=K.
\end{align}
\end{subequations}

To solve \eqref{SPQ}, we first solve \eqref{get_qu} in a closed-form by  first-order optimality condition:
{\small 
\begin{subequations}\label{q_tilde}
\begin{align}
\tilde{q}_u&=\Bigg[\sum\limits_{(I_i, I_j)\in \mathbb{B}} \left(\hat{p}_i-\hat{p}_j\right)\left(\hat{p}_i-\hat{p}_j\right)^{\top} \Bigg]^{-1}
\Bigg[ \sum_{(I_i, I_j)\in \mathbb{B}} (\hat{p}_i-\hat{p}_j)\left(a^*_{uij}-\left(\hat{b}_i-\hat{b}_j\right)\right)\Bigg]\\
&=\Bigg[\sum\limits_{(I_i, I_j)\in \mathbb{B}} \left(\hat{p}_i-\hat{p}_j\right)\left(\hat{p}_i-\hat{p}_j\right)^{\top} \Bigg]^{-1}
\Bigg[\sum_{(I_i, I_j)\in \mathbb{B}}
(\hat{p}_i-\hat{p}_j)\left(\left(\hat{p}_i-\hat{p}_j\right)^{\top}\hat{q}_u+(\varepsilon_{ui}-\varepsilon_{uj})\right)
\Bigg]\\
&=\hat{q}_u+\left(P_{\mathbb{B}}P_{\mathbb{B}}^{\top}\right)^{-1}P_{\mathbb{B}}\varepsilon_{u\mathbb{B}}, 
\end{align}
\end{subequations}}
where $P_{\mathbb{B}}$ is the matrix whose columns correspond to the latent factor vectors $\hat{p}_i-\hat{p}_j$ for all $(I_i, I_j)\in \mathbb{B}$, $\varepsilon_{u\mathbb{B}}$ is the vector whose elements correspond to the differences of noises $\varepsilon_{ui}-\varepsilon_{uj}$ for all $(I_i, I_j)\in \mathbb{B}$.
We assume that the noise term $\varepsilon_{ui}, i=1,\ldots,n$, follows i.i.d. Gaussian distribution ${N}(0, \sigma^2)$. 
Then we can write 
\begin{equation*}
    \tilde{a}_{uij}-\hat{a}_{uij}=(\hat{p}_{i}-\hat{p}_{j})^{\top}(\tilde{q}_u-\hat{q}_u)+\delta,\ i,j=1,\ldots,n, 
\end{equation*}
where $\delta\sim {N}(0, 4\sigma^2)$. 
Thus, by \eqref{q_tilde}, for any $u=1,\ldots,m,\ i,j=1,\ldots,n$, 
\begin{align*}
\mathbb{E}(\tilde{a}_{uij}-\hat{a}_{uij})^2&=\mathbb{E}\left[(\hat{p}_{i}-\hat{p}_{j})^{\top}(\tilde{q}_u-\hat{q}_u)\right]^2+\mathbb{E}[\delta^2]\\
&=\mathbb{E}\left[(\hat{p}_{i}-\hat{p}_{j})^{\top}\left(P_{\mathbb{B}}P_{\mathbb{B}}^{\top}\right)^{-1}P_{\mathbb{B}}\varepsilon_{u\mathbb{B}}\right]^2+
4\sigma^2, 
\end{align*}
\begin{align*}
   &\ \ \mathbb{E}\bigg[ \sum_{(I_i, I_j)\in \mathbb{I}^n\times \mathbb{I}^n}\left( \tilde{a}_{uij}-\hat{a}_{uij}\right)^2\bigg]
   =\mathbb{E}\left[\left\| P^{\top}\left(P_{\mathbb{B}}P_{\mathbb{B}}^{\top}\right)^{-1}P_{\mathbb{B}}\varepsilon_{u\mathbb{B}}\right\|_2^2\right]+
   4n^2\sigma^2 \\
   &  \overset{(1)}{=}\left\| P^{\top}\left(P_{\mathbb{B}}P_{\mathbb{B}}^{\top}\right)^{-1}P_{\mathbb{B}}C_{u\mathbb{B}}\right\|_F^2+
   4n^2\sigma^2 \overset{(2)}{=}2n^2\sigma^2\text{tr}\left(\left(P_{\mathbb{B}}P_{\mathbb{B}}^{\top}\right)^{-1}\right)+
   4n^2\sigma^2, 
\end{align*}
where $P$ is the matrix whose columns correspond to the latent factor vectors $\hat{p}_i-\hat{p}_j$ for all $(I_i, I_j)\in \mathbb{I}^n\times \mathbb{I}^n$ and $C_{u\mathbb{B}}$ denotes the square root of the covariance matrix of $\varepsilon_{u\mathbb{B}}$; 
{\color{black}equality $\overset{(1)}{=}$} follows from the definition of the Frobenius norm $\left\| \cdot \right\|_F$; 
{\color{black}equality 
$\overset{(2)}{=}$} holds {\color{black}by}  
$C_{u\mathbb{B}}=\sqrt{2\sigma^2}{E}_{K\times K}$ 
and 
{\color{black}the assumption that 
$PP^{\top}=n^2 
{E}_{v\times v}$, which reflects the isotropic position of the user's latent factor \cite{Anava2015}}.
Here, ${E}$ represents the identity matrix.

Then we can reformulate \eqref{SPQ} as
\begin{subequations}\label{re_SPQ}
\begin{align}
\min\limits_{\mathbb{B}}\ & \text{tr}\left[\left(P_{\mathbb{B}}P_{\mathbb{B}}^{\top}
\right)^{-1}
\right]\\
\text{s.t.}\ \ & \mathbb{B}\subset \mathbb{I}^n\times\mathbb{I}^n,\ \left|\mathbb{B}\right|=K.
\end{align}
\end{subequations}
By using a classical greedy iterative algorithm, for instance \cite{glover2013advanced,SPQ2018}, we can solve problem \eqref{re_SPQ}. 
The algorithm provides us with $K$ pairs of items, which form the questions in the questionnaire. 
Correspondingly, we denote the first item of the $k$-th pair as $W_k$, and the second item as $Y_k$. 

\section{Preference elicitation via pairwise comparison}\label{sec-pre-elic}
In this section, we introduce a novel optimization-based preference elicitation approach, taking advantage of the user's choices between pairwise 
items in {\color{black}the
questionnaire}. 
From the choices of the questions from the user, we can detect the user's preferences by comparing her/his expected utility of the items.
This is a major extension of the traditional pairwise comparison approach or choice-based model in the linear utility framework to the general non-parametric expected utility functions.

Suppose that we have generated a set of $K$ pairwise comparison questions through the SPQ algorithm, denoted as $\{(W_k,Y_k),\ k=1,\ldots,K\}$.  
Let 
$\mathscr{S} := \{0\}\cup\bigcup\limits_{k=1}^{K} \operatorname{supp}(W_k)\cup\bigcup\limits_{k=1}^{K} \operatorname{supp}(Y_k)\cup\{\bar{b}\} $, 
where $\bar{b}$ is the maximal outcome that the user {\color{black}cares}. 
We denote $N:=|\mathscr{S}|$ as the cardinality of set $\mathscr{S}$ and $\mathcal{N}:=[0, \bar{b}]$ as the domain of definition.  
Let $\overline{\mathbb{Y}} := \{\bar{y}_{j}\}_{j=1,\ldots,N}$ be the ordered sequence of points in $\mathscr{S}$ with fixed $\bar{y}_{1}=0$ and $\bar{y}_{N}=\bar{b}$.
We provide the user with 
{\color{black}the questionnaire consisting of $K$ questions}
and collect her/his choices.
For each 
{\color{black}question}
$(W_k,Y_k)$, $k=1,\ldots,K$, we denote $Z_k=1$ if the user chooses $W_k$, $Z_k=-1$ if she/he chooses $Y_k$, and $Z_k=0$ if she/he doesn't make a choice between $W_k$ and $Y_k$. 
We assume that all the choices are made by the user's true utility without any acquisition error.
Different from traditional conjoint analysis in which the user's utility is assumed to be linear to the features, we consider a non-parametric expected utility function $\mathbb{E}[u(\cdot)]$. 
\begin{assumption}\label{assump-1}
    The user's preference can be represented by Von Neumann--Morgenstern expected utility theory. 
\end{assumption}

We then aspire to estimate a reference utility function of the user according to her/his choices on the $K$ pairwise comparisons.
Let $\mathcal{L}^p({\mathcal{N}})$ denote
the set of real functions $u: \mathcal{N}\to \mathbb{R}$ integrable to the $p$-th order.
Denote ${L}_c$ as a subset of $\mathcal{L}^p({\mathcal{N}})$ consisting of all twice continuous differential, monotone increasing, concave and normalized functions, i.e., 
$${L}_c:= \{ u\in\mathcal{L}^2({\mathcal{N}})\mid u'(y)\geq 0, u''(y)\leq 0, u(0)=0, u(\bar{b})=1,\ \forall y\in\mathcal{N} \}.$$ 
{\color{black}It directly follows from Assumption \ref{assump-1} that the true utility function of the user is contained in set ${L}_c$.}
Then we can introduce a set 
$$L_{K} := \{ u\in {L}_c \mid Z_k \cdot \mathbb{E}\left[u(W_k) \right] \geq Z_k \cdot \mathbb{E}\left[u(Y_k) \right],\ k = 1, \ldots, K \},$$ 
representing all utility functions in ${L}_c$ consistent 
{\color{black}with}
the choices of the user.
{\color{black}The constraints in $L_{K}$ indicate that if the user prefers $W_k$ to $Y_k$, then $\mathbb{E}\left[u(W_k)\right] \geq \mathbb{E}\left[u(Y_k)\right]$, and if the user prefers $Y_k$ to $W_k$, then $\mathbb{E}\left[u(W_k)\right] \leq \mathbb{E}\left[u(Y_k)\right]$.}
Since we have assumed that the choices are made without any 
acquisition error, the true utility of the user ought to be in the set $L_{K}$.
However, due to the insufficient number of 
{\color{black}questions}, $L_{K}$ is commonly not a one-point set, i.e., the utility functions consistent with the user's choices are not unique. 
{\color{black}
\begin{lemma}
    The set $L_{K}$ is a one-point set if and only if 
\begin{equation}\label{lemma:L_K}
\sup\limits_{u\in L_k}\mathbb{E}[u(\xi)]=\inf\limits_{u\in L_k}\mathbb{E}[u(\xi)],\ \forall\xi\in \mathcal{L}^p(\Omega, \mathcal{F}, P).    
\end{equation}
\end{lemma}
\begin{proof}
If $L_{K}$ is a one-point set, then the equation holds trivially. Conversely, suppose the equality \eqref{lemma:L_K} holds for every $\xi\in \mathcal{L}^p(\Omega, \mathcal{F}, P)$. Now assume, to obtain a contradiction, that there exist $u_1,u_2\in L_{K}$ such that $u_1\neq u_2$, i.e., there exists $y'\in \mathcal{N}$ such that $u_1(y')\neq u_2(y')$. Let $\xi^{'}\in \mathcal{L}^p(\Omega, \mathcal{F}, P)$ be a particular random variable taking the value of $y'$ with probability of 1. Consequently, $\mathbb{E}[u_1(\xi^{'})]\neq\mathbb{E}[u_2(\xi^{'})]$, leading to $\sup_{u\in L_k}\mathbb{E}[u(\xi^{'})]\neq\inf_{u\in L_k}\mathbb{E}[u(\xi^{'})]$, which contradicts the equality \eqref{lemma:L_K}.
\end{proof}
}

\begin{figure}[htbp]
	\centering 
    \includegraphics[scale=0.5]{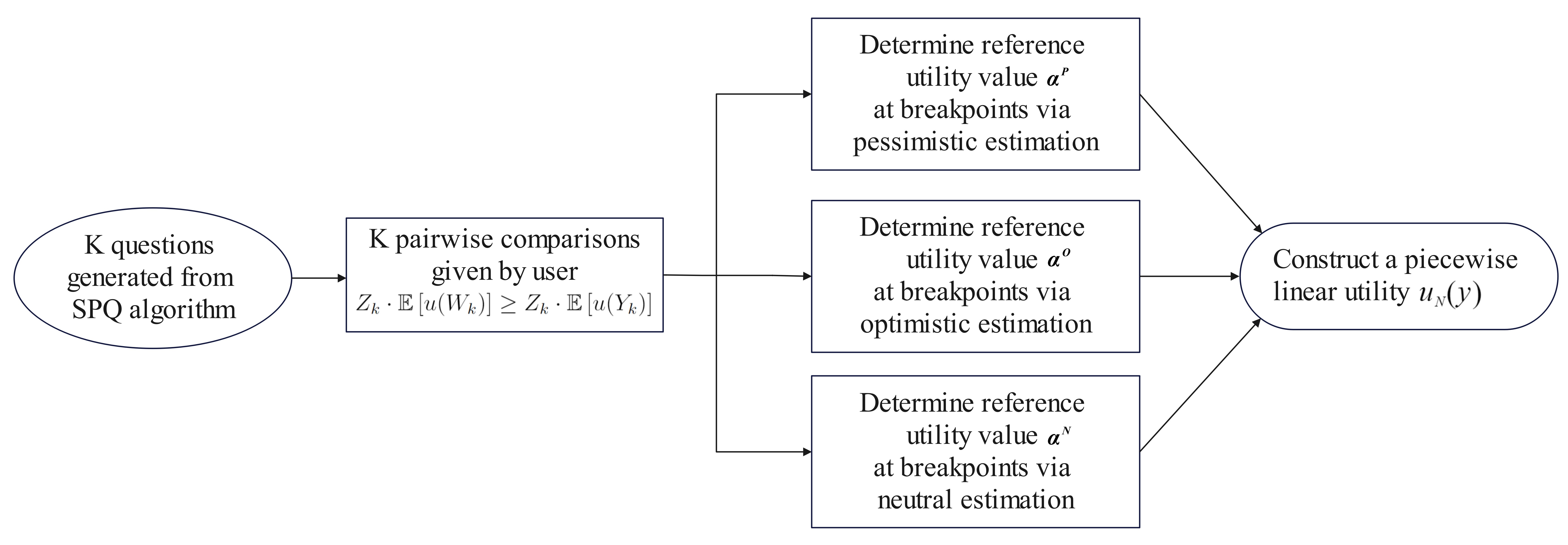}
        \caption{{\color{black}Nominal preference elicitation process}}
        \label{PE-flowchart}
\end{figure}

Thus, we shall put forward some criteria to select one reference utility function from the set $L_{K}$ for portfolio recommendation, which we refer to as a nominal utility function. 
In subsections \ref{sec-pess}--\ref{sec-unbi}, we introduce three kinds of estimation approaches to determine the nominal utility function, namely the pessimistic estimation, the optimistic estimation, and the {\color{black}neutral} estimation. 
The procedures are summarized in a flowchart in Figure~\ref{PE-flowchart}. 

\subsection{Pessimistic estimation}\label{sec-pess}
As stated above, commonly the utility functions consistent 
{\color{black}with}
the user's choices are not unique, we now design some criteria to select one nominal utility function from $L_{K}$.
The first approach is to find a pessimistic estimation by choosing the worst-case utility function in $L_{K}$ with the smallest expected reference utility value. 

To do so, we first specify a random benchmark $\tilde{y}$ and compute the worst-case utility over $L_K$ by minimizing the expected utility of $\tilde{y}$, i.e., we construct 
{\color{black}the}
following preference elicitation optimization problem: 
\begin{equation*}
    \inf \limits_{u\in L_k} \mathbb{E}[u(\tilde{y})]. 
\end{equation*}
For practical use, we {\color{black} can select}
the benchmark $\tilde{y}$ based on items in the 
{\color{black}questionnaire}. 
Let $\eta :=\left[I_1, I_2, \ldots, I_n\right]^{\top}$ denote the random vector where the $i$-th component corresponds to 
item $I_i$ in the item set $\mathbb{I}^n$. 
Let $\tilde{x}\in \mathbb{R}^n$ be a pre-determined constant vector that can be viewed as a {\color{black}benchmark} portfolio on the item set. 
We then set $\tilde{y}:=h(\tilde{x}, \eta)=\tilde{x}^{\top}\eta$ as the random benchmark and use the following optimization problem to determine a pessimistic estimation of the user's nominal utility: 
\begin{equation}\label{PRO}
    \inf\limits_{u\in L_K}\mathbb{E}\left[u(h(\tilde{x}, \eta)) \right].   
\end{equation}
We assume that all the components in $\eta$ are independent of each other, thus the joint distribution of $\eta$ is the product of all the marginal distribution functions of $I_i, i=1,\ldots,n$. 
Because of the finite support of all the items $I_i$, the random vector $\eta$ has finite scenarios. 
We assume that there are $M$ scenarios of $\eta$, denoted by $\Omega := \{\eta_1,\ldots,\eta_M\}$ 
with associated probabilities $p_{i} := \mathbb{P}\left[\eta=\eta_i\right]$. 
Thus, $\mathbb{E}\left[u(h(\tilde{x}, \eta))\right]=\sum_{i=1}^{M}p_i\cdot u(h(\tilde{x},\eta_i))$. 

However, problem \eqref{PRO} requires optimization over an infinite dimensional decision variable space $L_K$, which is not easy to address from a computational perspective. 
To reduce the dimension of the decision variable, we utilize the piecewise linear approximation method, i.e., we approximate $L_K$ by a subset $L_N\subset L_K$, where $L_N$ consists of some piecewise linear functions $u_N(\cdot)$ defined over the interval $\left[\bar{y}_1, \bar{y}_N\right]$ with breakpoints on $\overline{\mathbb{Y}}$. 
As the number of pieces increases, the difference between $L_N$ and $L_K$ progressively reduces to zero \cite{MS2021}. 
Let $\alpha:=\left[\alpha_1,\ldots,\alpha_N\right]\in\mathbb{R}^N$ denote the utility values at the breakpoints in $\overline{\mathbb{Y}}$
and $\beta:=\left[\beta_1,\ldots,\beta_{N-1}\right]\in\mathbb{R}^{N-1}$ denote the slopes of linear pieces between two adjacent breakpoints, i.e., $\alpha_j:=u_N(\bar{y}_j),\ j=1,\ldots, N$ and $\beta_j:=(\alpha_{j+1}-\alpha_j)/(\bar{y}_{j+1}-\bar{y}_j),\ j=1,\ldots,N-1$. 
Then any function $u_N$ in $L_N$ can be structured as 
\begin{equation}\label{utility function}
    u_N(y)=\left\{ 
    \begin{aligned}
    &1 && y= \bar{y}_N,\\
    &\beta_j(y-\bar{y}_j)+\alpha_j && \bar{y}_j\leq y <\bar{y}_{j+1},\ j=1,\ldots,N-1,\\
    & 0 && y=\bar{y}_1.\\
    \end{aligned}
    \right.
\end{equation}
By limiting utility functions as piecewise linear, we {\color{black}approximately} solve problem \eqref{PRO} by 
\begin{equation}\label{approximated PRO}
\inf\limits_{u_N\in L_N}\sum\limits_{i=1}^M p_i \cdot u_N(h(\tilde{x},{\color{black}\eta_i})). 
\end{equation}

Since the function $u_N\in L_N \subset L_K$ is required to be monotone increasing and concave, we should have $\beta_{j+1}\leq\beta_j,\ j=1,\ldots, N-2$, and $\beta_j\geq0,\ j=1,\ldots,N-1$. 
Additionally, due to the piecewise linear structure of $u_N$, we have the following equivalent reformulation of $u_N$ at any point $y\in[\bar{y}_1, \bar{y}_N]$: 
\begin{subequations}\label{u_N:reformulation_1}
\begin{align}
u_N(y)=\min\limits_{v\geq0, w\in\mathbb{R}}\
& vy+w   \\
 \text{s.t.}\ \ & v\bar{y}_j+w\geq\alpha_j,\ j=1,\ldots,N.
\end{align}
\end{subequations}
By applying 
reformulation \eqref{u_N:reformulation_1} {\color{black}to}
the objective function and adding the constraints in $L_N$, we can reformulate problem \eqref{approximated PRO} as the following 
linear programming problem: 
\begin{subequations}\label{inf_linear}
\begin{align}
\label{obj}
\min\limits_{\alpha, \beta, v, w}\
& \sum_{i=1}^M p_i(v_ih(\tilde{x},\eta_i)+w_i)   \\ 
 \text{s.t.}\ \
\label{constraints_inf}
& v_i\bar{y}_j+w_i\geq\alpha_j,\ i=1,\ldots,M,\ j=1,\ldots,N, \\ 
\label{constraints_noml}
& \alpha_1 =0, \alpha_N=1,\\
\label{constraints_user}
& Z_k \cdot \sum\limits_{j=1}^N \mathbb{P}\left[W_k=\bar{y}_j\right]\alpha_j \geq Z_k \cdot \sum\limits_{j=1}^N \mathbb{P}\left[Y_k=\bar{y}_j\right]\alpha_j,\ k=1,\ldots,K, \\ 
\label{constraints_alpha_beta}
& \alpha_{j+1}=\beta_j(\bar{y}_{j+1}-\bar{y}_j)+\alpha_j,\ j=1,\ldots,N-1, \\
\label{constraints_beta}
& \beta_{j+1}\leq\beta_j,\ j=1,\ldots,N-2, \\
\label{constraints_domain}
& v\in \mathbb{R}_+^M,\ w\in\mathbb{R}^M,\ \alpha\in\mathbb{R}^N,\ \beta\in\mathbb{R}_+^{N-1}. 
\end{align}
\end{subequations}
In \eqref{inf_linear}, \eqref{constraints_noml} is the normalization constraint, \eqref{constraints_alpha_beta} reflects the relation between $\alpha$ and $\beta$, \eqref{constraints_beta} characterizes the concavity of the function. 
Non-negativeness of $\beta$ in \eqref{constraints_domain} guarantees the monotonicity of the optimal utility function. 
\eqref{constraints_user} requires the optimal utility function 
{\color{black}to}
be consistent with the user's choices over {\color{black}$K$ questions in the questionnaire}. 
We denote the optimal solution of \eqref{inf_linear} as $\alpha^P=\left[\alpha_1^P,\ldots,\alpha_N^P\right]$ and $\beta^P=\left[\beta_1^P,\ldots,\beta_{N-1}^P\right]$, which form a pessimistic nominal utility function $u^P(\cdot)$. 

\subsection{Optimistic estimation}\label{sec-opti}
The second approach is to seek an optimistic estimation by choosing a utility function consistent with the user's pairwise {\color{black}comparison} choices with the 
{\color{black}highest}
expected reference utility value. 
Building upon the parameters and notations introduced in subsection \ref{sec-pess}, we formulate the {\color{black}optimistic estimation problem} as 
\begin{equation}\label{sup}
    \sup\limits_{{\color{black} u_N}\in L_N}\mathbb{E}\left[{\color{black}u_N}(h(\tilde{x}, \eta)) \right]. 
\end{equation}

Different from \eqref{u_N:reformulation_1}, we apply another reformulation of $u_N\in L_N$ {\color{black} as the minimum of $N-1$ linear functions}
\begin{equation}\label{u_N:reformulation_2}
    u_N(y)=\min\limits_{j=1,\ldots,N-1}\left[\beta_j(y-\bar{y}_j)+\alpha_j\right], 
\end{equation}
which holds due to the piecewise linear and concave property of $u_N$. 
By applying \eqref{u_N:reformulation_2} and taking account all the constraints in $L_N$, we reformulate \eqref{sup} as the following finite dimensional linear programming problem:
\begin{subequations}\label{sup_linear}
\begin{align}
\max\limits_{\alpha,\beta,z}\
& \sum_{i=1}^M p_i z_i   \\
 \text{s.t.}\ \  & z_i\leq\beta_j(h(\tilde{x},\eta_i)-\bar{y}_j)+\alpha_j,\ i=1,\ldots,M, \ j=1,\ldots,N-1, \\
\label{constraint_start}
& \alpha_1 =0,\ \alpha_N=1,\\ 
& Z_k \cdot \sum\limits_{j=1}^N\mathbb{P}\left[W_k=\bar{y}_j\right]\alpha_j \geq Z_k \cdot \sum\limits_{j=1}^N\mathbb{P}\left[Y_k=\bar{y}_j\right]\alpha_j,\ k=1,\ldots,K, \\
& \alpha_{j+1}=\beta_j(\bar{y}_{j+1}-\bar{y}_j)+\alpha_j,\ 
j=1,\ldots,N-1, \\
& \beta_{j+1}\leq\beta_j,\ j=1,\ldots,N-2, \\
\label{constraints_end}
& \alpha\in\mathbb{R}^N,\ \beta\in\mathbb{R}_+^{N-1},\ z\in\mathbb{R}^M.  
\end{align}
\end{subequations}
Here the constraints \eqref{constraint_start}--\eqref{constraints_end} reflect the constraint $u_N\in L_N$, {\color{black} same as \eqref{constraints_noml}--\eqref{constraints_domain}}. 
$z=[z_1, z_2, \ldots, z_M]$ are auxiliary variables which characterize the minimal value 
{\color{black}among all $N-1$ linear functions.}
We denote $\alpha^O=\left[\alpha_1^O,\ldots,\alpha_N^O\right]$ and $\beta^O=\left[\beta_1^O,\ldots,\beta_{N-1}^O\right]$ as the optimal solution of \eqref{sup_linear}, which form an optimistic nominal utility function $u^O(\cdot)$.

\subsection{Neutral estimation}\label{sec-unbi}
The third approach aims to find a utility function in $L_N$ that lies 
{\color{black}at the midpoint between}
the pessimistic nominal utility function and the optimistic nominal utility function in the sense of a certain distance. 
As the distances from the 
{\color{black} third}
estimation to the optimistic 
{\color{black} one}
and the pessimistic 
{\color{black} one}
are the same, we name it the neutral estimation. 

We start by defining a semi-distance metric to measure the distance between two utility functions. 
Let $\mathscr{G}$ be a set of measurable functions defined over the interval $\left[\bar{y}_1,\bar{y}_N\right]$.
For any $u,v\in L_K$, we define the semi-distance between $u$ and $v$ as
$$\dd_{\mathscr{G}}(u,v):=\sup_{g\in \mathscr{G}}\left| \int_{\bar{y}_1}^{\bar{y}_N}g(z)du(z) - \int_{\bar{y}_1}^{\bar{y}_N}g(z)dv(z)\right|,$$ 
where the integrals are taken in the sense of Lebesgue-Stieltjes integration. 
Here, $g$ can be thought of as a test function and $\dd_{\mathscr{G}}(u, v)=0$ means that for all test functions in $\mathscr{G}$, there is no difference between $u$ and $v$ albeit that $u\neq v$. 
Since in $L_K$ we normalize the utility functions with $\alpha_1=0$,\ $\alpha_N=1$, the semi-distance $\dd_{\mathscr{G}}(u,v)$ takes the form of the pseudo-metric on the $\zeta$-structure in probability theory. 
As a special case, we focus on the Kantorovich metric, denoted by $\dd_K(u,v)$, by taking 
\fontsize{11}{15}{
\begin{equation*}
\mathscr{G}=\mathscr{G}_L:=
    \{g:\left[\bar{y}_1,\bar{y}_N\right]\rightarrow [0, 1] \mid g \text{ is Lipschitz continuous with modulus bounded by 1} \}. 
\end{equation*}
}
Since $u,v$ are in normalized utility set $L_c$, we have
\begin{equation*}
\int_{\bar{y}_1}^{\bar{y}_N} g(z) d u(z)-\int_{\bar{y}_1}^{\bar{y}_N} g(z) d v(z)=\int_{\bar{y}_1}^{\bar{y}_N}(g(z)-g(\bar{y}_1)) d u(z)-\int_{\bar{y}_1}^{\bar{y}_N}(g(z)-g({\bar{y}_1})) d v(z).
\end{equation*}
We then can replace $\mathscr{G}_L$ with
\fontsize{10}{15}{
\begin{equation}\label{eq:kantor-GL}
\mathscr{G}_L:=
    \{g:{\color{black}\left[{\bar{y}_1},\bar{y}_N\right]}\rightarrow \mathbb{R} \mid g({\bar{y}_1})=0,\ g \text{ is Lipschitz continuous with modulus bounded by 1} \}. 
\end{equation}
}
In the forthcoming discussions, we consider $\mathscr{G}_L$ defined in \eqref{eq:kantor-GL}.

We then find a utility function $u$ which lies 
{\color{black}at the midpoint between}
the pessimistic nominal utility function $u^P$ and optimistic nominal utility function $u^O$ by minimizing the larger distance from $u$ to $u^P$ or $u^O${\color{black}:}
\begin{equation}\label{avg_minmax}
    \min\limits_{u\in L_K}\max\{\dd_K(u,{u}^{P}),\dd_K(u,{u}^{O})\}.
\end{equation}
{\color{black}If}
\eqref{avg_minmax} reaches its optimal solution $u^*$, $\dd_K(u^*, u^P)=\dd_K(u^*, u^O)$ and $u^*$ locates at the center of $u^P$ and $u^O$. 

As both $u^P$ and $u^O$ are piecewise linear functions with breakpoints in $\bar{\mathbb{Y}}$, the optimal $u^*$ keeps the same piecewise linear structure~\cite{MS2021}. 
Thus, we only consider $u$ in $L_N$. 
The integral of $g$ over $u$ in $L_N$ can be calculated as 
\begin{equation*}
\int_{\bar{y}_1}^{\bar{y}_N}g(z)du(z)=\sum\limits_{j=1}^{N-1}\beta_j
\int_{\bar{y}_j}^{\bar{y}_{j+1}}g(t)dt.
\end{equation*}
And for any $u,\tilde{u}\in L_N$, 
\begin{equation*}
\dd_K(u,\tilde{u})=\sup\limits_{g\in\mathscr{G}_L}\sum\limits_{j=1}^{N-1}
    (\beta_j-\tilde{\beta}_j)\int_{\bar{y}_j}^{\bar{y}_{j+1}}g(t)dt,
\end{equation*}
where $\beta_j$ ($\tilde{\beta_j}$, resp.) is the slope of $u$ ($\tilde{u}$, resp.) in the interval $\left[\bar{y}_j,\bar{y}_{j+1}\right]$, $j=1,\ldots,N-1$. 

{\color{black}
By \cite[Proposition 6]{MS2021}, we have that the Kantorovich distance between two piecewise linear {\color{black}functions} in $L_N$ can be computed by the following quadratically constrained quadratic program (QCQP): 


\begin{subequations}
\label{eq:Kant-u-v-LP-new-1}
\begin{align}
\dd_K(u, \tilde{u})=\max\limits_{\substack{w_1, \ldots, w_{N-1}\\ z_2, \ldots, z_N}}\ &\sum_{j=1}^{N-1}  (\beta_j-\tilde{\beta}_j)w_j\\
\text{s.t.}\quad &  
 w_{j}\leq \frac{1}{4} (\bar{y}_{j+1}- \bar{y}_{j})^2 -\frac{1}{4}( z_{j+1}-z_{j})^2 + \frac{1}{2}(\bar{y}_{j+1}- \bar{y}_{j})( z_{j+1}+z_{j}),\\
&-w_{j}\leq \frac{1}{4}(\bar{y}_{j+1}- \bar{y}_{j})^2 -\frac{1}{4}( z_{j+1}-z_{j})^2 - \frac{1}{2}(\bar{y}_{j+1}- \bar{y}_{j})( z_{j+1}+z_{j}),\\
&\text{for}  \ j=1,\ldots,N-1, \nonumber
\end{align}
\end{subequations}
where $z_1=g(\bar{y}_1)=0$.
\noindent Considering that problem \eqref{eq:Kant-u-v-LP-new-1} is convex but not strictly convex with respect to $z_j$, we make a variable transformation to $x_{j}=z_{j+1}-z_j$, $j=1,\ldots,N-1$ and subsequently formulate \eqref{eq:Kant-u-v-LP-new-1} as
\begin{subequations}
\label{eq:Kant-u-v-LP-new-x}
\begin{align}
\dd_K(u, \tilde{u})=\max\limits_{\substack{w_1, \ldots, w_{N-1}\\ x_1, \ldots, x_{N-1}}}\ &\sum_{j=1}^{N-1}  (\beta_j-\tilde{\beta}_j)w_j\\
\text{s.t.}\quad  &  
 w_{j}\leq \frac{1}{4} (\bar{y}_{j+1}- \bar{y}_{j})^2 -\frac{1}{4}( x_{j})^2 + \frac{1}{2}(\bar{y}_{j+1}- \bar{y}_{j})(  2\sum_{i=1}^{j-1}x_i+x_j ),\\
&-w_{j}\leq \frac{1}{4}(\bar{y}_{j+1}- \bar{y}_{j})^2 -\frac{1}{4}( x_{j})^2 - \frac{1}{2}(\bar{y}_{j+1}- \bar{y}_{j})(2\sum_{i=1}^{j-1}x_i+x_j ),\\
&\text{for}  \ j=1,\ldots,N-1.\nonumber
\end{align}
\end{subequations}
Since \eqref{eq:Kant-u-v-LP-new-x} is a strictly convex quadratic programming problem, it can be represented as a self-dual second-order cone program:
\begin{subequations}
\label{eq:Kant-u-v-LP-new-socp}
\begin{align}
\shortintertext{$\dd_K(u, \tilde{u})=$}
\max\limits_{\substack{w_1, \ldots, w_{N-1}\\ x_1, \ldots, x_{N-1}}}\ &\sum_{j=1}^{N-1}  (\beta_j-\tilde{\beta}_j)w_j\\
\text{s.t.}\  &  
\left[\begin{array}{c}
     x_j \\
     -2w_{j} + (\bar{y}_{j+1}- \bar{y}_{j})( 2\sum_{i=1}^{j-1}x_i+x_j ) + \frac{1}{2}(\bar{y}_{j+1}- \bar{y}_{j})^2  - \frac{1}{2}  \\
     -2w_{j} + (\bar{y}_{j+1}- \bar{y}_{j})(2\sum_{i=1}^{j-1}x_i+x_j ) + \frac{1}{2}(\bar{y}_{j+1}- \bar{y}_{j})^2 + \frac{1}{2} 
\end{array}\right] \in \text{SOC}_{3},\\ 
  &  
\left[\begin{array}{c}
     x_j \\
     2w_{j} - (\bar{y}_{j+1}- \bar{y}_{j})(2\sum_{i=1}^{j-1}x_i+x_j ) + \frac{1}{2}(\bar{y}_{j+1}- \bar{y}_{j})^2  - \frac{1}{2}  \\
     2w_{j} - (\bar{y}_{j+1}- \bar{y}_{j})(2\sum_{i=1}^{j-1}x_i+x_j ) + \frac{1}{2}(\bar{y}_{j+1}- \bar{y}_{j})^2 + \frac{1}{2} 
\end{array}\right] \in \text{SOC}_{3}, \\ 
&\text{for}  \ j=1,\ldots,N-1.\nonumber
\end{align}
\end{subequations}
where $\text{SOC}_{3}:=\{v\in\mathbb{R}^3\mid\|[v_1,v_2]^{\top}\|\leq v_3\}$ denotes the 3-dimensional second-order cone. 
Taking duality gives 
\begin{subequations}\label{kant:dual-socp}
\fontsize{10}{15}{
\begin{align}
\shortintertext{$\dd_K(u, \tilde{u})=$}
\min_{t,\lambda,\mu,s,\rho,\phi} 
\ &\frac{1}{2}\sum_{j=1}^{N-1}\left(\bar{y}_{j+1}-\bar{y}_{j}\right)^2 \left(\lambda_j+\mu_j+\rho_j+\phi_j\right)-\frac{1}{2} \sum_{j=1}^{N-1}\left(\lambda_j-\mu_j+\rho_j-\phi_j\right) \\
\text { s.t. }\ &
\label{kant:dual-const-1}
\frac{1}{2}\left(\beta_j-\tilde{\beta}_j\right) - \lambda_j-\mu_j+\rho_j+\phi_j=0, \\
\label{kant:dual-const-2}
& t_j+s_j+\sum_{i=j+1}^{N-1} 2\left(\bar{y}_{i+1}-\bar{y}_i\right)\left(\lambda_i+\mu_i-\rho_i-\phi_i\right)+\left(\bar{y}_{j+1}-\bar{y}_j\right)\left(\lambda_j+\mu_j-\rho_j-\phi_j\right)=0, \\
&\left[t_j, \lambda_j, \mu_j\right]^{\top} \in \text{SOC}_3,\left[s_j, \rho_j, \phi_j\right]^{\top} \in \text{SOC}_3, \\
&\text{for }\ j=1, \ldots, N-1.\nonumber
\end{align}}
\end{subequations}
}
As we consider $u\in L_N$, we have the following constraints for the slopes $\beta$ and the utility values $\alpha$ at breakpoints:
\begin{subequations}\label{constraint_L_N}
    \begin{align}
     &\beta_j=(\alpha_{j+1}-\alpha_j)/(\bar{y}_{j+1}-\bar{y}_j),\ j=1,\ldots,N-1, \\
     &\beta_{j+1} \leq \beta_j,\ j=1,\ldots,N-2,  \\
     &\alpha_1=0,\ \alpha_N=1,\\
     &\alpha\in\mathbb{R}^N,\ \beta\in\mathbb{R}^{N-1}_+.
     \end{align}
\end{subequations}
{\color{black}
Employing the dual form of $\dd_K(u,\tilde{u})$ in \eqref{kant:dual-socp} and the constraints in \eqref{constraint_L_N}, we can reformulate 
problem \eqref{avg_minmax} as the following 
second-order cone programming problem: 
{\fontsize{10}{15}\selectfont
\begin{subequations}\label{final avg}
\begin{align}
\min\limits_{\substack{\alpha,\beta,\zeta,t,s,\\ \lambda,\mu,\rho,\phi}}\
& \zeta   \\
 \text{s.t.}\ \  
& \zeta \geq \frac{1}{2}\sum_{j=1}^{N-1}\left(\bar{y}_{j+1}-\bar{y}_{j}\right)^2 \left(\lambda_j^k+\mu_j^k+\rho_j^k+\phi_j^k\right)-\frac{1}{2} \sum_{j=1}^{N-1}\left(\lambda_j^k-\mu_j^k+\rho_j^k-\phi_j^k\right),\ k=1,2,\\
& \frac{1}{2}\left(\beta_j-\beta^P_j\right) - \lambda_j^1-\mu_j^1+\rho_j^1+\phi_j^1=0,\ j=1,\ldots,N-1,\\
& \frac{1}{2}\left(\beta_j-\beta^O_j\right) - \lambda_j^2-\mu_j^2+\rho_j^2+\phi_j^2=0,\ j=1,\ldots,N-1,\\
& t_j^k+s_j^k+\sum_{i=j+1}^{N-1} 2\left(\bar{y}_{i+1}-\bar{y}_i\right)\left(\lambda_i^k+\mu_i^k-\rho_i^k-\phi_i^k\right)+\left(\bar{y}_{j+1}-\bar{y}_j\right)\left(\lambda_j^k+\mu_j^k-\rho_j^k-\phi_j^k\right)=0,\nonumber\\
& j=1,\ldots,N-1,\ k=1,2,\\
&\left[t_j^k, \lambda_j^k, \mu_j^k\right]^{\top} \in \text{SOC}_3,\left[s_j^k, \rho_j^k, \phi_j^k\right]^{\top} \in \text{SOC}_3, \\
& Z_k \cdot \sum\limits_{j=1}^N\mathbb{P}\left[W_k=\bar{y}_j\right]\alpha_j \geq Z_k \cdot \sum\limits_{j=1}^N\mathbb{P}\left[Y_k=\bar{y}_j\right]\alpha_j,\ k=1,\ldots,K, \\
&(\bar{y}_{j+1}-\bar{y}_j)\cdot\beta_j=(\alpha_{j+1}-\alpha_j),\ j=1,\ldots,N-1 \\
&\beta_{j+1} \leq \beta_j,\ j=1,\ldots,N-2,  \\
&\alpha_1=0,\ \alpha_N=1,\\
& \alpha\in\mathbb{R}^N,\ \beta\in\mathbb{R}^{N-1}_+,\ \zeta\in\mathbb{R}_+,\ t^k,s^k,\lambda^k,\mu^k,\rho^k,\phi^k\in\mathbb{R}^{N-1},\ k=1,2.
\end{align}
\end{subequations}}
}
We denote the optimal solution of \eqref{final avg} as $\alpha^N=[\alpha_1^N,\ldots,\alpha_N^N]$ and $\beta^N=[\beta_1^N,\ldots,\beta_{N-1}^N]$, which form the neutral nominal utility function $u^N(\cdot)$ in the shape of \eqref{utility function}. 

\section{Utility maximization}\label{sec-portfolio-opt}
Once the robo-advisor has obtained the user's personalized nominal utility function $u^X(\cdot)$, where $X$ could be one of the {\color{black}estimations} $P$, $O$, $N$ derived through the preference elicitation process in Section \ref{sec-pre-elic},  the robo-advisor {\color{black}recommends} an investment portfolio according to the user's nominal {\color{black}preference}. 
Suppose there are $S$ risky assets with random return rates $\xi_1, \xi_2,\ldots, \xi_S$ and 
{\color{black}a risk-free asset with a deterministic return rate $\xi_0=0$.} 
{\color{black}Denote $\xi=[\xi_0,\xi_1,\ldots,\xi_S]
$.} 
The robo-advisor determines a portfolio $x=[x_0,x_1,x_2,\ldots,x_S]\in\mathbb{R}^{S+1}$ by maximizing the user's expected nominal utility of the portfolio return $\xi^{\top}x$ with a total investment budget $W_0$ of the user: 
{\color{black}\begin{subequations}\label{port_org}
\begin{align}
\mathop{\max}\limits_{x\in \mathbb{R}^{S+1}}\
& \mathbb{E}\left[u^X(\xi^{\top} x) \right]   \\ 
 \text{s.t.}\ \  & \sum_{s=0}^S x_s=W_0,\\
\label{port_org:const-c}
& 0\leq x_s\leq c_s ,\ s=0, \ldots, S, 
\end{align}
\end{subequations}}
where $c_s$ is the maximal investment limit on the $s$-th asset. 

Suppose that we have collected $T$ historical return rate data of the $S$ risky assets, {\color{black}which can be viewed as $T$ samples of the random return $\xi$.}
Let $\xi_t:= \left[{\color{black}\xi_0,}\xi_1^t, \xi_2^t, \ldots, \xi_S^t\right]$ where $\xi_s^t$ is the return rate of the $s$-th risky asset in the $t$-th sample. 
By using the historical samples, we can solve a sample average approximation problem of \eqref{port_org}: 
{\color{black}
\begin{subequations}\label{port_avg}
\begin{align}
\mathop{\max}\limits_{x\in \mathbb{R}^{S+1}}\
& \frac{1}{T}\sum_{t=1}^T u^X({\xi_t}^\top x)
\\ 
\text{s.t.}\ \  & \sum_{s=0}^S x_s=W_0,\\
& 0\leq x_s\leq c_s ,\ s=0, \ldots, S. 
\end{align}
\end{subequations}}
Noticing that $u^X,\ X=P,O,N$, obtained in Section \ref{sec-pre-elic} is a {\color{black}piecewise linear function} with slopes $\beta_j^X$ and function values $\alpha_j^X$ at the {\color{black}left end} point of the $j$-th piece, $j=1,\ldots,N-1$.
Note that $u^X$ is concave, thus, 
\begin{equation*}
u^X(y)=\min\limits_{j=1,\ldots,N-1} \beta_j^X(y-\bar{y}_j)+\alpha_j^X,\ X=P,O,N.
\end{equation*}
By introducing auxiliary variables $z=[z_1, z_2, \ldots, z_T]^\top$, we reformulate 
problem \eqref{port_avg} as a linear programming problem: 
{\color{black}
\begin{subequations}\label{port_final}
\begin{align}
\max\limits_{x\in \mathbb{R}^{S+1}{\color{black},z\in\mathbb{R}^T}}\
& \frac{1}{T}\sum_{t=1}^{T} z_t \\
 \text{s.t.}\ \qquad  & z_t\leq \beta_j^X(\xi_t^{\top} x-\bar{y}_j)+\alpha_j^X,\ t=1,\ldots, T,\ j=1,\ldots, N-1,\\
& \sum_{s=0}^S x_s=W_0,\\
& 0\leq x_s\leq c_s,\ s=0, \ldots, S.
\end{align}
\end{subequations}
}
The robo-advisor then 
{\color{black}offers}
the optimal solution of problem \eqref{port_final} to the user as the personalized investment strategy. 

\section{Numerical test}\label{sec-numerical}
In this section, we carry out a series of numerical tests to evaluate the performance of the proposed robo-advisor system. 
The tests consist of two parts: simulation experiments and 
{\color{black} human experiments.}
The simulation experiments assume a closed-form true utility function and 
1) compare the efficiency of static preference questionnaires (SPQ) and randomly generated questionnaires in subsection \ref{sec-sim-spq}; 
2) validate the convergence of nominal utility function to the true utility function in subsection \ref{sec-sim-utility}; 
3) simulate a long-term investment process 
in subsection \ref{sec-sim-invest}. 
For {\color{black} human experiments,}
{\color{black}
we have invited
three {\color{black}human} participants in subsection \ref{sec-real-utility} and 60 {\color{black}human} participants in subsection \ref{sec-real-port} to experience the robo-advisor {\color{black} system} for verifying practical feasibility.}

\subsection{Test setting}\label{sec-test-sett}
When carrying out simulation experiments, we generate pairwise comparisons from a set of 20 items denoted by $\mathbb{I}^{20}$ with a maximum outcome of \text{\textyen}500,000 (i.e., $\bar{y}_N=500000$), 
{\color{black}
as shown in Table~\ref{20-item set} in supplementary materials. }
{\color{black} In} the questionnaire generation process (SPQ), we use some quadratic functions of the form $f(x)=ax^2+bx+c$ to simulate the historical ratings required as the input of the SPQ. 
Here, we generate $a, b, c$ from uniform {\color{black} distributions} within the interval [-50, 50]. 
For each randomly generated $\tilde f(x)=\tilde{a}x^2+\tilde{b}x+\tilde{c}$, we compute $\mathbb{E}\left[\tilde{f}(I_i)\right],\ i=1,\ldots,20$, and normalize the {\color{black} expected utility values} to a scale of 0 to 10, which are used as the {\color{black} historical} ratings of the 20 items in $\mathbb{I}^{20}$. 
In simulation experiments (subsection \ref{sec-sim-spq}--subsection \ref{sec-sim-invest}), we preset a virtual user with 
{\color{black}the}
true utility $u^T(y):=\frac{1}{1-e^{-5}}(1-e^{-10^{-5}y})$, which is increasing, concave and normalized to [0, 1] within the interval [0, 500000]. 

{\color{black}
When implementing applications to {\color{black} human} participants (subsection \ref{sec-sim-utility}--subsection \ref{sec-real-port}),
}
we generate pairwise comparisons from a smaller set of 10 items $\mathbb{I}^{10}$ to reduce the interaction difficulties, as shown in Table \ref{10-item set}. 
We collect ratings for the 10 items in $\mathbb{I}^{10}$ from random {\color{black} human} individuals, with scores ranging from 0 to 10, which are used in the questionnaire generation process (SPQ). 
Moreover, it's worth mentioning that 
{\color{black}
the robo-advisor system 
does not collect
the participants' personal information 
during the experiment, and the participants are unaware of the underlying principles of the robo-advisor system. 
}

\subsection{Simulation results: comparison between SPQ and randomly generated questionnaires}\label{sec-sim-spq}

Based on the item set $\mathbb{I}^{20}$ shown in Table \ref{20-item set}, we compare the efficiency of the SPQ method with purely randomly generated questionnaires, by examining the fitness of the nominal utility functions elicited by the same elicitation models with different questionnaires generated by SPQ or randomly. 
{\color{black}In detail, we draw pair of items from the item set $\mathbb{I}^{20}$ for $K$ times to create the randomly generated questionnaire.
The SPQ method is discussed in Section \ref{sec-spq}.
The virtual user $u^T(\cdot)$ then fills out these two questionnaires with $K$ pairwise comparison questions based on her/his expected utility comparison. Then we can elicit three nominal utility functions accordingly from the answers.}

\begin{figure}[htbp]
    \centering
    \includegraphics[scale=0.4]{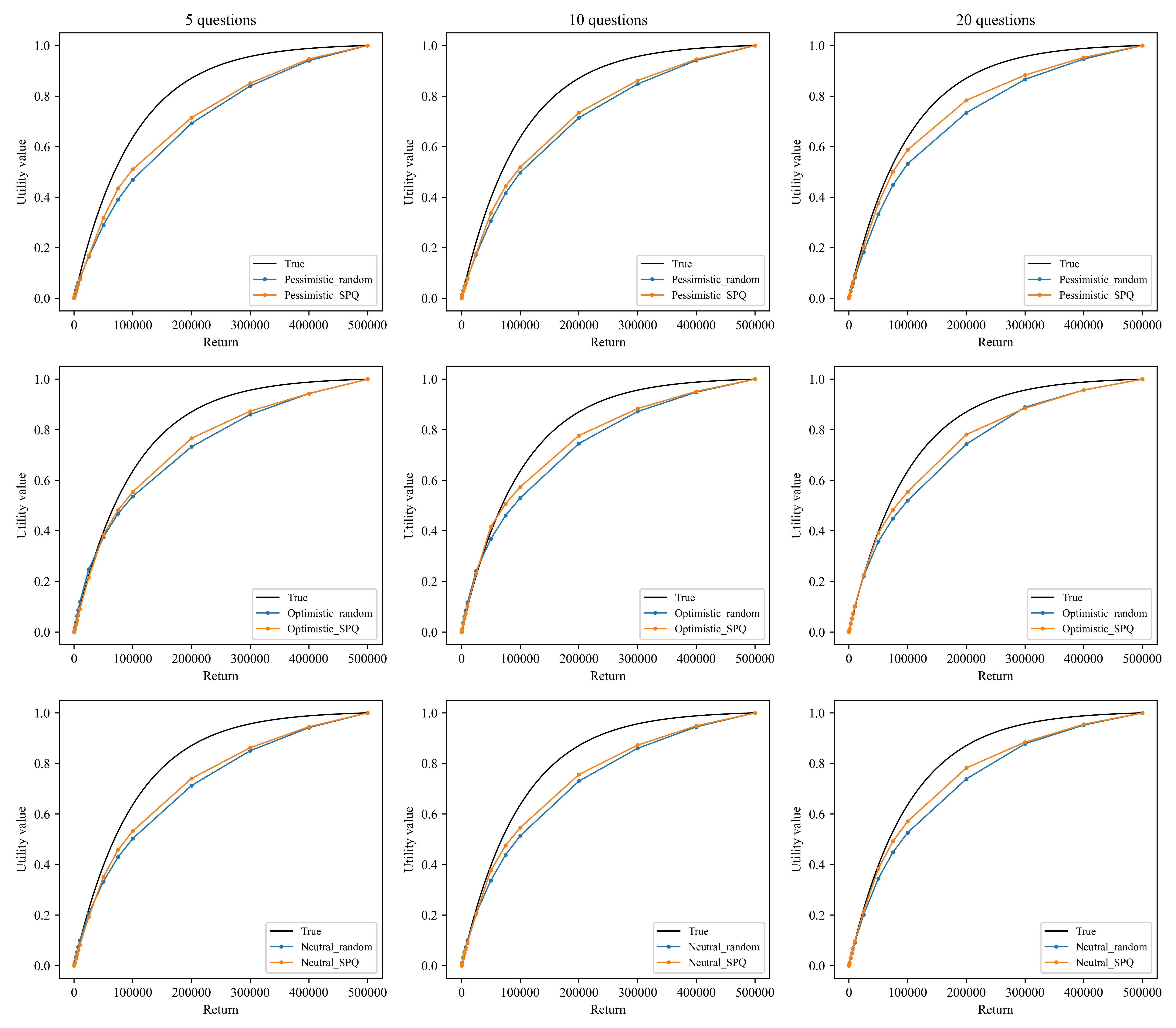}
    \caption{Nominal utility functions (pessimistic, optimistic, neutral estimations) elicited by using SPQ and randomly generated questionnaires}
    \label{fig:comparison_spq}
\end{figure}

We set a group of case studies with different numbers of 
{\color{black}questions}
$K=5, 10, 20$ and different elicitation models. 
{\color{black}The nominal/true utility functions are presented in Figure \ref{fig:comparison_spq} and the Kantorovich distances between these nominal utility functions and true utility functions are provided in Table \ref{tab:kantor}.}
To eliminate randomness in {\color{black} questionnaire} generation, all preference elicitation 
in Figure~\ref{fig:comparison_spq} and Table \ref{tab:kantor} are derived from {\color{black}repeat} independent experiments {\color{black} for} 50 times and taking the average.

\begin{table}[htbp]
  \centering
  \caption{\color{black}Kantorovich distances between nominal utility functions and true utility function}
  {\fontsize{8}{12}\selectfont\color{black}
    \begin{tabular}{cc|ccc}
    \toprule
          &       & 5 questions & 10 questions & 20 questions \\
    \midrule
    \multirow{2}[2]{*}{Pessimistic-True} & SPQ   & 0.0941 & 0.0855 & 0.0552 \\
          & Random & 0.1108 & 0.0989 & 0.0824 \\
    \midrule
    \multirow{2}[2]{*}{Optimistic-True} & SPQ   & 0.0665 & 0.0564 & 0.0578 \\
          & Random & 0.0799 & 0.0749 & 0.0724 \\
    \midrule
    \multirow{2}[2]{*}{Neutral-True} & SPQ   & 0.0803 & 0.0700 & 0.0563 \\
          & Random & 0.0943 & 0.0862 & 0.0773 \\
    \bottomrule
    \end{tabular}}
  \label{tab:kantor}%
\end{table}%


By Figure~\ref{fig:comparison_spq} and Table \ref{tab:kantor}, we can make the following observations:

\begin{itemize}[itemsep=2pt,topsep=0pt,parsep=0pt]
    \item As the number of questions in the questionnaire increases, the nominal utility functions derived from both SPQ and randomly generated questionnaires fit the true utility function $u^T$ better under any of the three preference elicitation models (pessimistic, optimistic, neutral). 
    \item With the same number of questions, the nominal utility functions derived from SPQ fit the true utility function better, compared to those obtained from randomly generated questionnaires. 
    \item 
    In cases with fewer questions, SPQ demonstrates a significant advantage over randomly generated questionnaires in preference elicitation. 
\end{itemize}

\subsection{Simulation results: convergence of nominal utility functions to 
{\color{black}the}
true utility function}
\label{sec-sim-utility}
\begin{figure}[htbp]
    \centering
    \includegraphics[scale=0.45]{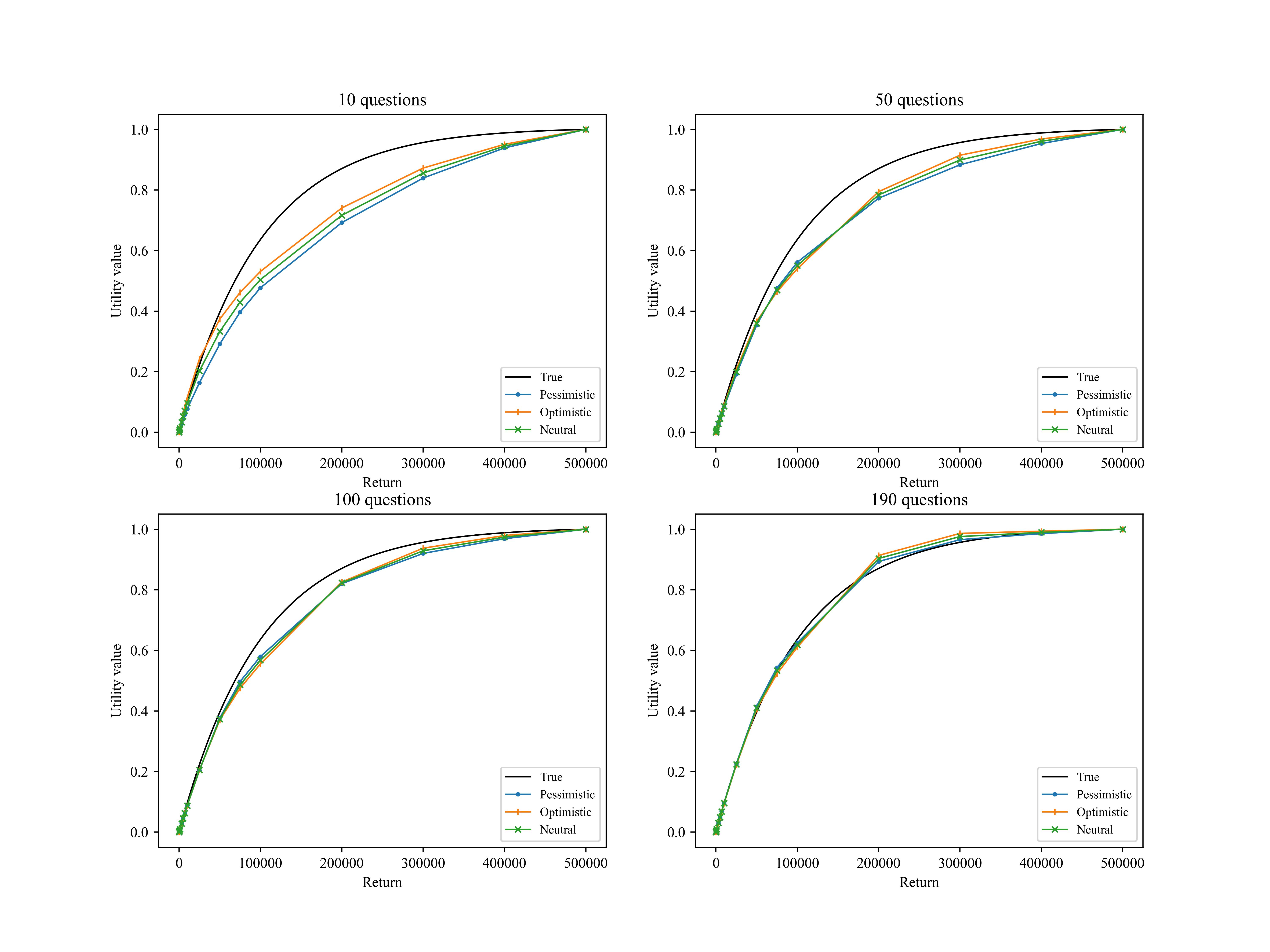}
    \caption{Nominal utility functions (pessimistic, optimistic, neutral estimation) elicited by different numbers of questions generated randomly}
    \label{fig:fitting}
\end{figure}

\begin{figure}[htbp]
    \centering
    \includegraphics[scale=0.6]{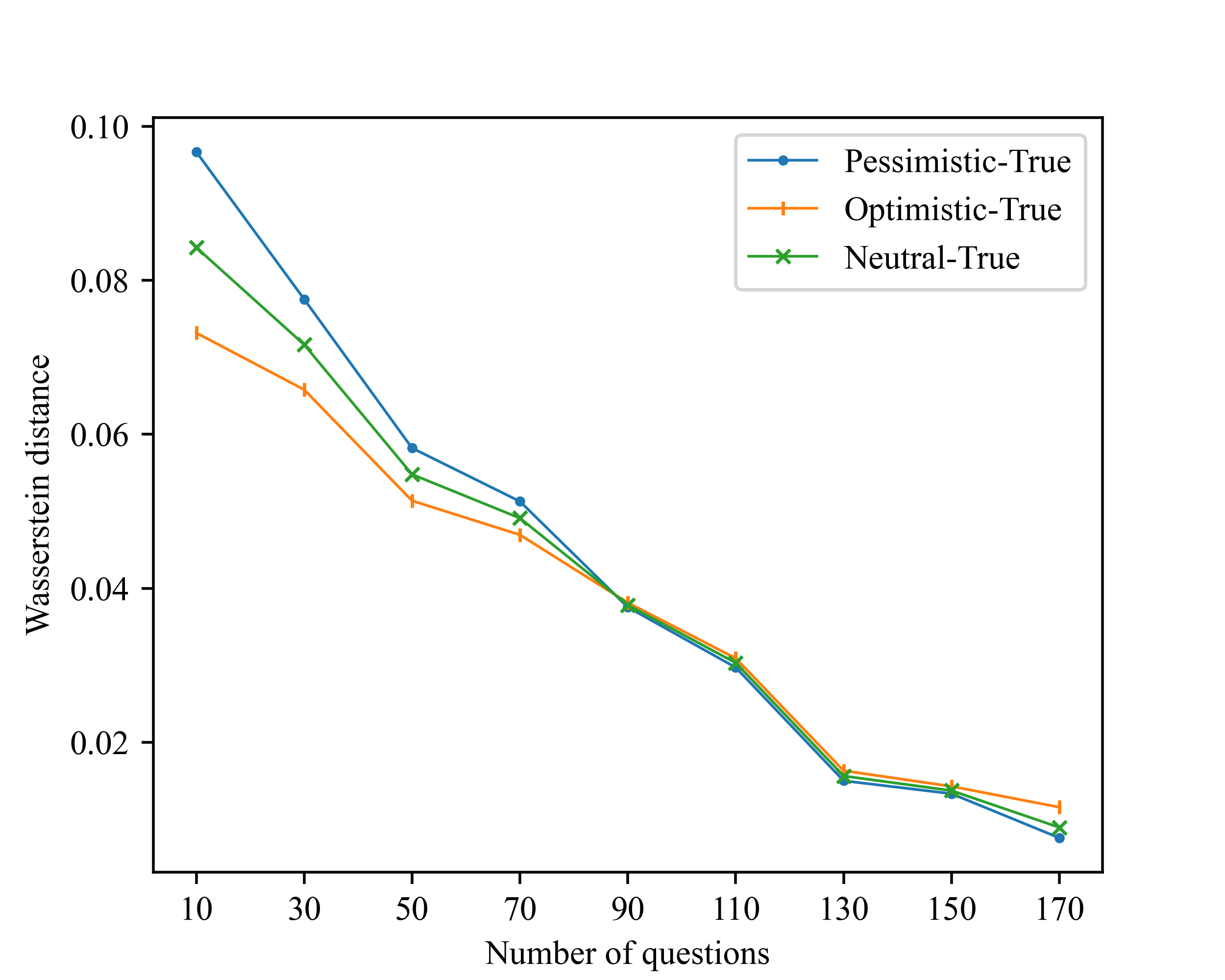}
    \caption{Kantorovich distance between nominal utility functions and the true utility function when the number of questions increases}
    \label{fig:wasserstein_global}
\end{figure}

To examine the efficiency of the proposed elicitation models, we investigate whether the nominal utility functions $u^X(\cdot), X=P,O,N$, generated by the three preference elicitation models converge to the preset true utility function $u^T$. 
We {\color{black} study the randomly generated questionnaires only and }
extend the numerical test with 
more randomly generated questions. 
Figure~\ref{fig:fitting} shows the fitness of the nominal utility functions estimated from $K$=10, 50, 100 and 190 questions by using pessimistic, optimistic and neutral estimation methods. 
Figure~\ref{fig:wasserstein_global} shows the Kantorovich distance between nominal utility functions and true utility function with increasing number of questions. 
Similarly, to reduce the impact of randomness in {\color{black} questionnaire} generation, all the preference elicitation results in Figures~\ref{fig:fitting},
\ref{fig:wasserstein_global} 
are derived from repeating independent experiments {\color{black} for} 100 times and taking the average.

We can find {\color{black} from Figure \ref{fig:fitting} and Figure \ref{fig:wasserstein_global}} the following conclusions: 
\begin{itemize}[itemsep=2pt,topsep=0pt,parsep=0pt]
    \item The curves of the nominal utility functions elicited by the neutral estimation method typically lie between those elicited by the pessimistic estimation method and the optimistic estimation method, which is consistent with the theoretical results. 
    \item As the number of questions increases, all three nominal utility functions (pessimistic, optimistic, neutral) gradually converge to the true utility function. 
\end{itemize}

\subsection{Simulation results: investment process driven by the robo-advisor}\label{sec-sim-invest}
Upon obtaining an estimation of the nominal utility, the robo-advisor proceeds to compute a portfolio for recommendation by solving the linear programming problem~\eqref{port_final}. 
In this subsection, we select a set of funds and examine the investment performance of the portfolios recommended by the robo-advisor in a long-term investment process, compared with the portfolio determined by the true utility. 


{\color{black}In this case study, we consider an asset universe with 5 exchange-traded-funds (ETF) 
and one risk-free asset.
In detail, the first two ETFs are representative of the energy (XLE) and the technology (XLK) sectors within the S\&P500 market index (SPY). These two 
are based on the industry partition of SPY, which is adopted as the benchmark in the optimization problem. Other three ETFs, which are poorly or negatively correlated with the benchmark, are also considered to facilitate portfolio diversification: the SPDRGold Shares (GLD), which tracks the performance of gold bullion, an ETF for long-term (7 to 10 years) treasury bond investments (IEF), and finally an ETF constructed to track the US dollar performance (WisdomTree Bloomberg U.S. Dollar Bullish Fund, USDU). Additionally, we have a cash account with null return. We assume no transaction costs, and the results are generated with constants $c_s=0.4\cdot W_0$, $s=1,\ldots S$ in the specification of problem \eqref{port_final}, where $W_0$ is the pre-set initial investment budget. 
We collect daily returns of the five ETFs from January 3, 2017 to December 29, 2023\footnote{The data are downloaded from https://cn.investing.com/.}. 
We conducted Adjusted Dickey-Fuller (ADF) tests on the daily returns of the five ETFs during this period, all demonstrating time series stability. 
}

\begin{figure}[htbp]
    \centering
    \includegraphics[scale=0.7]{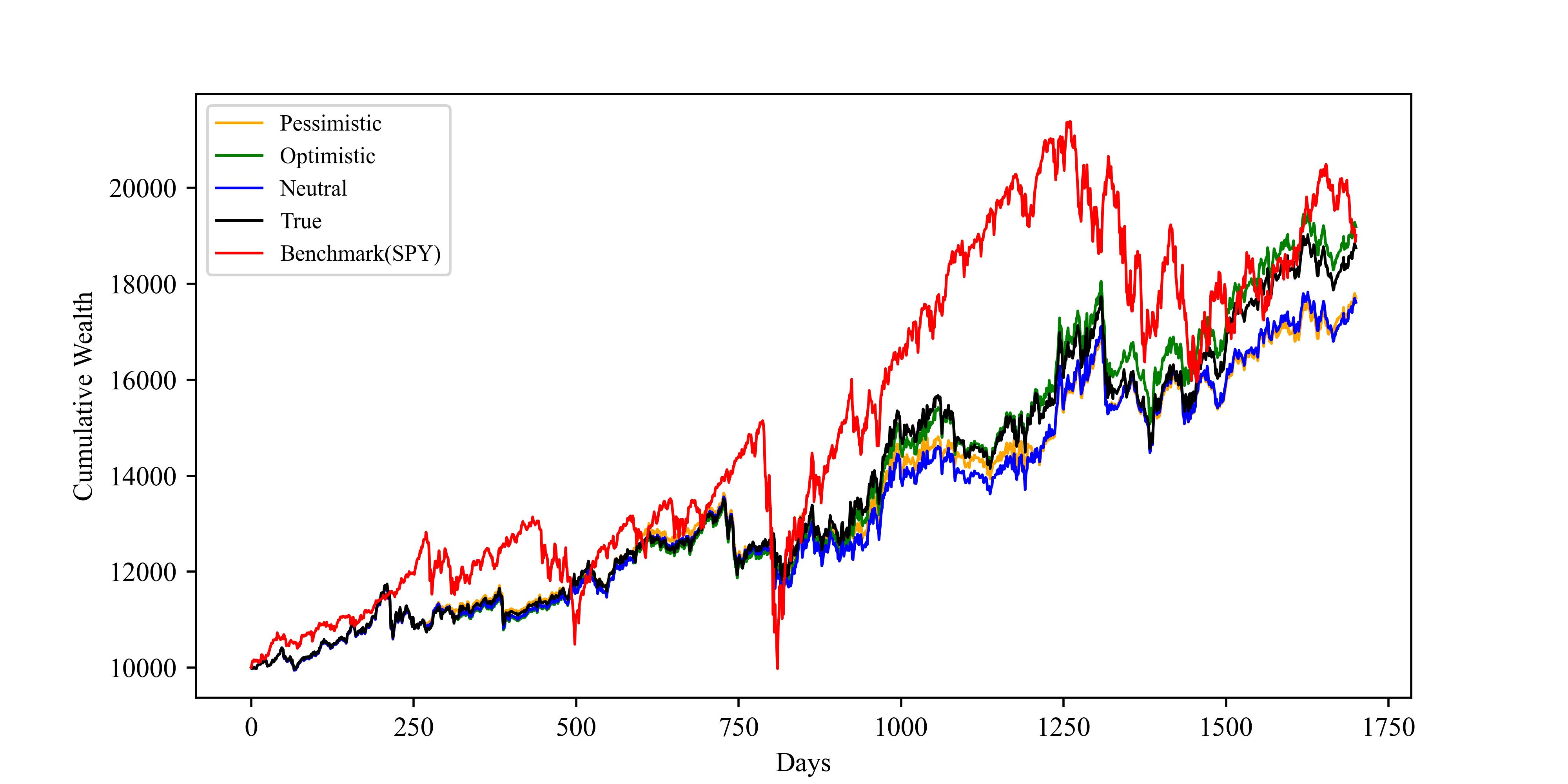}
    \caption{{\color{black}Out-of-sample cumulative wealth by investing with nominal utility functions (pessimistic, optimistic, neutral) and the true utility function }}
    \label{fig:portfolio}
\end{figure}

{\color{black}We utilize the static preference questionnaire which selects 15 pairs of items from the Item set $\mathbb{I}^{20}$ shown in Table \ref{20-item set}} to elicit three nominal utility functions. The user's choice data and the elicitation methods are the same as those 
{\color{black}in subsection \ref{sec-sim-spq}.}
We set the initial investment budget as 10,000 (i.e., $W_0=10000$). 
The investment process is simulated through a rolling time window methodology, 
{\color{black} and the portfolio is re-balanced for each seven days. }
{\color{black} In detail, }
we first set a historical in-sample date period of 
{\color{black}60 days from January 3, 2017 to March 29, 2017 (i.e., $T=60$),} 
then we compute three optimal portfolios by solving problem~\eqref{port_final} with the three 
nominal utility functions. 
Meanwhile, we compute a portfolio by solving \eqref{port_final} with the true utility function as a benchmark.  
Following that, the four portfolios are used to invest on the 
{\color{black}five ETFs}
on the subsequent seven days, i.e., we compute four out-of-sample returns based on the real daily returns 
{\color{black}from March 30, 2017 to April 7, 2017.}
We then roll the time window 
{\color{black}seven days forward,}
by removing 
{\color{black}the first seven days}
from the window and adding 
{\color{black}the next seven days}
to the window. 
Repeating the procedure above until reaching the last day in the dataset, 
{\color{black}December 29, 2023,}
we obtain four series of 
{\color{black} daily returns}
corresponding to the four nominal/true utility functions. 
We calculate the corresponding cumulative wealth curves and show them in Figure~\ref{fig:portfolio}. 

{\color{black}We can observe that the portfolios recommended by the robo-advisor are generally aligned with the true preferences of the user, as evidenced by the close resemblance between the cumulative wealth curves generated by the three nominal utility functions and the true utility function.
And, as expected, the generated portfolios diversify risks from the market itself with a similar final expected return but obviously smaller volatility, compared to the benchmark (SPY).
}


\subsection{Application to {\color{black} human} users: elicited utility}\label{sec-real-utility}
In this subsection, we apply the robo-advisor to three {\color{black} human} participants to examine the efficiency of the preference elicitation. 
The three representative participants are selected 
with a significant difference. 
We present some of their personal characteristics in Table~\ref{tab:user information}. 
{\color{black} The personal characteristics are not inputted to the robo-advisor system, but only used for post-analysis.}

\begin{table}[htbp]
\setlength{\extrarowheight}{1pt}
    \centering
    \caption{Personal characteristics of three participants}
    {\fontsize{10}{15}\selectfont
    \resizebox{\textwidth}{!}{
    \begin{tabular}{l l l l}
    \toprule
    &  User 1 & User 2 & User 3\\
    \midrule
    Gender & Male & Female & Male\\
    Age & 20 years old & 45 years old & 32 years old\\
    Occupation & Undergraduate student & \parbox[t]{5cm}{Production and Transportation \\ [-1.3ex] Equipment Operator} & \parbox[t]{3.5cm}{Business and Service \\ [-1.3ex] Industry Professional}\\
    Monthly Income & Below \text{\textyen}1,000 & \text{\textyen}5,000--\text{\textyen}10,000 & Above \text{\textyen}10,000\\
    Investment Experience & No & Little & Rich\\
    Investment Budget & \text{\textyen}5,000 & \text{\textyen}80,000 & \text{\textyen}1,000,000\\
    \bottomrule
    \end{tabular}}}
    \label{tab:user information}
\end{table}

\begin{figure}[h]
  \centering
  \subfloat[User 1]{\includegraphics[width=0.3\linewidth]{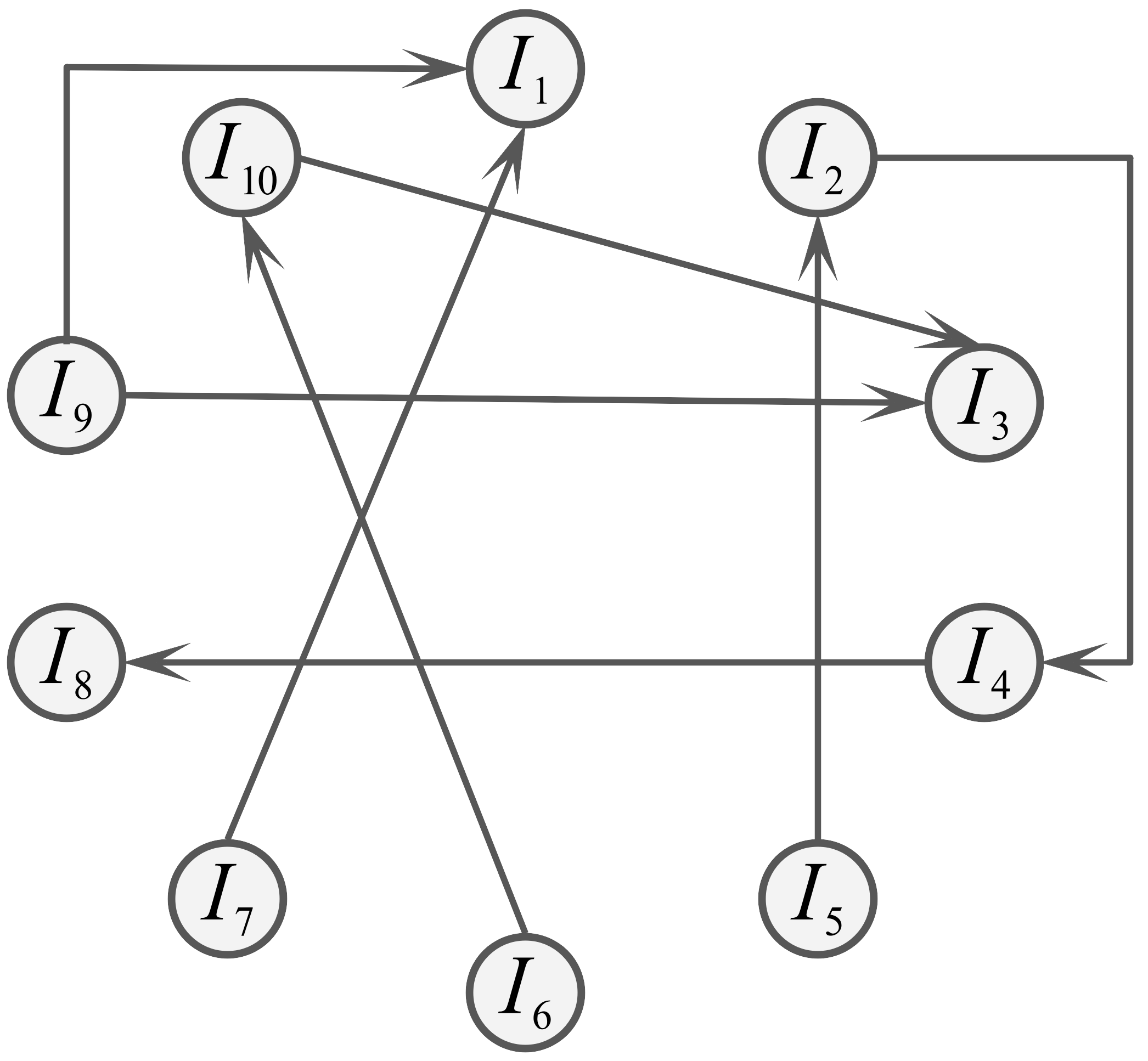}}\quad
  \subfloat[User 2]{\includegraphics[width=0.3\linewidth]{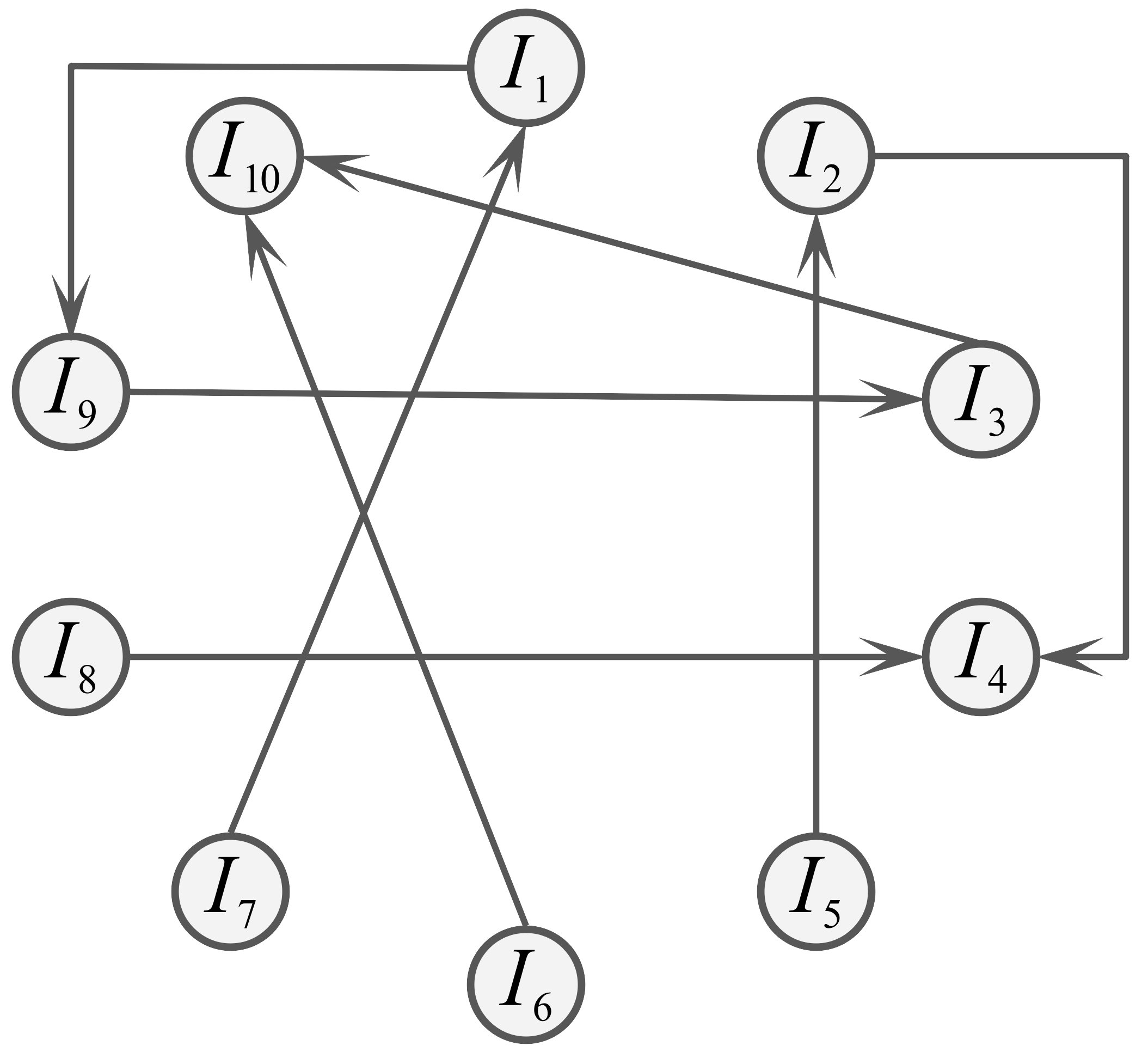}}\quad
  \subfloat[User 3]{\includegraphics[width=0.3\linewidth]{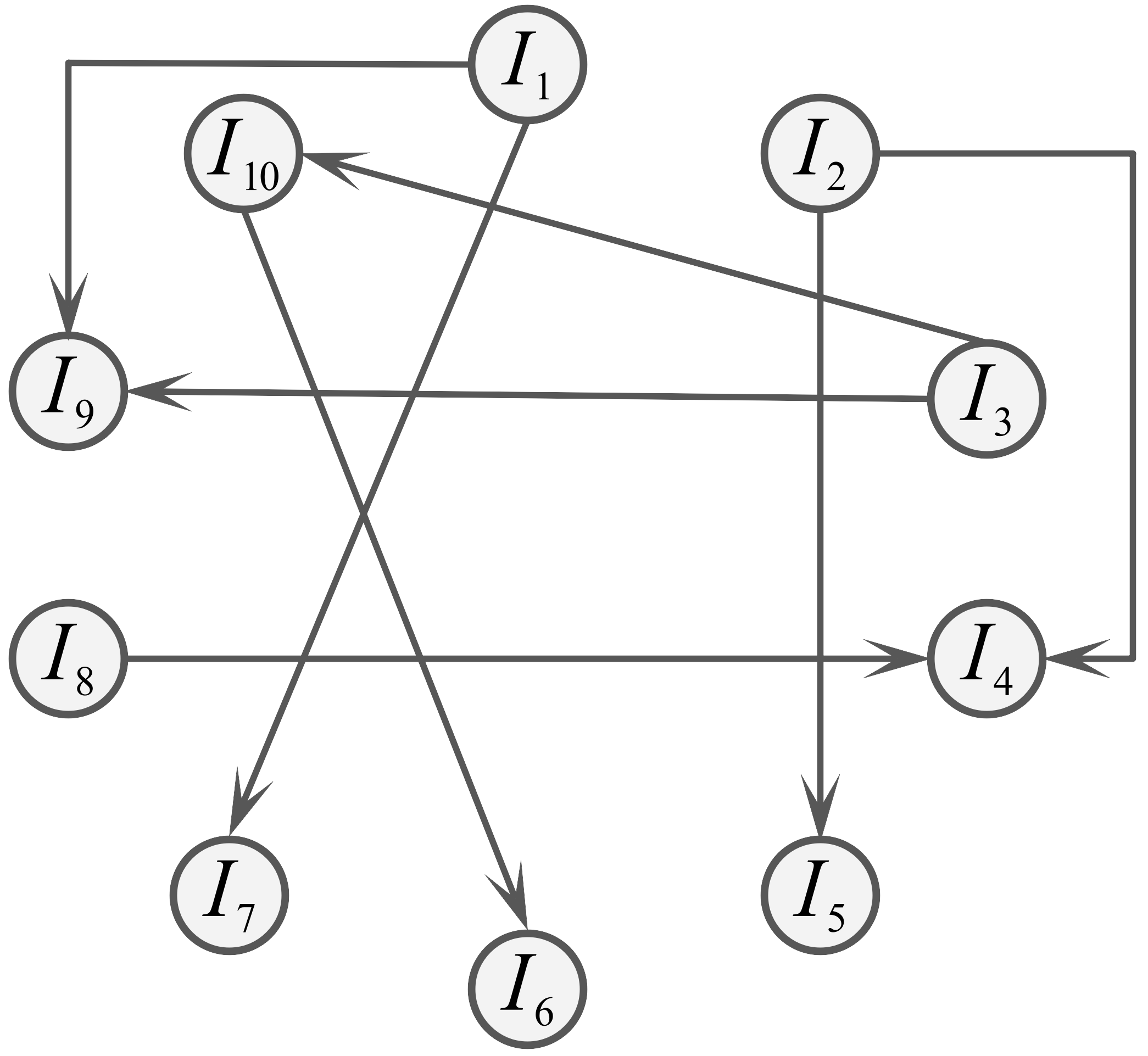}}
  \caption{\color{black}Preference graphs of three participants' answers to the preference questionnaire}
  \label{fig:graph user}
\end{figure}

As introduced in subsection \ref{sec-test-sett}, we use a 10-item set $\mathbb{I}^{10}$ and employ the SPQ method to generate a questionnaire consisting of 8 pairwise items: $(I_1, I_7)$, $(I_4, I_2)$, $(I_2, I_5)$, $(I_6, I_{10})$, $(I_3, I_9)$, $(I_4, I_8)$, $(I_1, I_9)$ and $(I_3, I_{10})$. 
We then ask the three participants to answer the questionnaire independently and record their choices for each pairwise items.  
{\color{black}The answers are provided in the three preference graphs respectively in Figure \ref{fig:graph user}. 
Here, a preference graph \cite{vcaklovic2017universal} is a tuple $\mathcal{G}=(V,A)$ where $V$ is a set of $n$ vertices (representing the items), $A \subseteq V \times V$ is a set of $m$ directed edges. Edge$(I_i, I_j)$ points from $I_i$ to $I_j$ if and only if the participant prefers $I_j$ over $I_i$. 
}

\begin{figure}[htbp]
    \centering
    \includegraphics[scale=0.4]{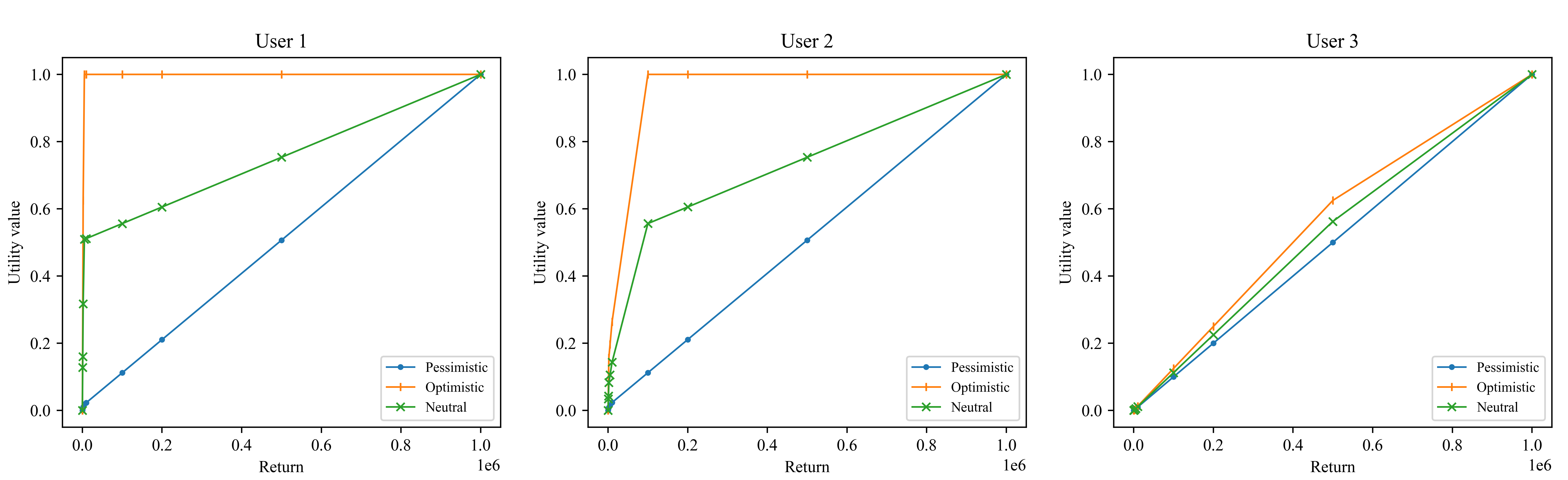}
    \caption{Nominal utility functions of the three participants under pessimistic, optimistic and neutral estimations}
    \label{fig:user}
\end{figure}

Based on the binary choices from the three participants, we construct three nominal utility functions for each participant through pessimistic, optimistic, and neutral estimation approaches proposed in Section \ref{sec-pre-elic}. 
We show the participants' nominal utility function curves under three estimation models in Figure~\ref{fig:user}. 

{\color{black}
Furthermore, according to the computational method for non-smooth utility function in \cite[Remark 3.9]{wurth2011risk}, we calculate the absolute risk aversion (ARA) and relative risk aversion (RRA) characteristics at breakpoints for the three users as
$${\rm RRA}(y)=\frac{y[u'(y^{-})-u'(y^+)]}{2u'(y^{-})},\quad 
{\rm ARA}(y)=\frac{u'(y^{-})-u'(y^+)}{2u'(y^{-})},\quad
\text{for } y \in \overline{\mathbb{Y}}\backslash\{\bar{y}_1, \bar{y}_N\},
$$
where $u'(y^{-})$ is the left-derivative of $u(\cdot)$ on $y$, $u'(y^+)$ is the right-derivative of $u(\cdot)$ on $y$.
These calculations of RRA and ARA characteristics are applied to the three users' neutral nominal utility functions at pre-set breakpoints $\overline{\mathbb{Y}}\backslash\{\bar{y}_1, \bar{y}_N\}=$\{800, 1000, 2000, 5000, 10000, 100000, 200000, 500000\}, 
shown in Table \ref{tab:user-RRA}. 

\begin{table}[htbp]
  \centering
  \caption{\color{black}RRA and ARA  characteristics at pre-set breakpoints of the neutral nominal utility functions for the three users}
  \resizebox{\textwidth}{!}{\color{black}
    \begin{tabular}{c|rrrrrrrr|rrrrrrrr}
    \hline
          & \multicolumn{7}{c}{RRA}   &       & \multicolumn{8}{c}{ARA} \\
          \hline 
          Breakpoints & 800 & 1000 & 2000 & 5000 & 10000 & 100000 & 200000 & 500000 & 800 & 1000 & 2000 & 5000 & 10000 & 100000 & 200000 & 500000 \\
    \hline
    User 1 & 312.71  & 20.81  & 0.00  & 58.08  & 4936.67  & 0.00  & 0.00  & 0.00  & 0.39  & 0.02  & 0.00  & 0.01  & 0.49  & 0.00  & 0.00  & 0.00  \\
    User 2 & 0.00  & 30.46  & 810.87  & 0.00  & 1978.89  & 44662.30  & 0.00  & 0.00  & 0.00  & 0.03  & 0.40  & 0.00  & 0.19  & 0.44  & 0.00  & 0.00  \\
    User 3 & 0.00  & 0.00  & 0.00  & 0.00  & 0.00  & 0.00  & 0.00  & 57522.12  & 0.00  & 0.00  & 0.00  & 0.00  & 0.00  & 0.00  & 0.00  & 0.11  \\
    \hline
    \end{tabular}}
  \label{tab:user-RRA}%
\end{table}%

}

We can observe from {\color{black}Figure \ref{fig:user} and Table \ref{tab:user-RRA}} that: 
The proposed robo-advisor model is highly applicable to {\color{black} human} users. 
The number of questions in the static preference questionnaire is acceptable and efficient for normal users. 
The robo-advisor estimates distinctive nominal utility functions 
{\color{black}and coefficients of risk aversion}
to different users, which are consistent with users' personal characteristics 
in Table~\ref{tab:user information}.
{\color{black} In detail,} we can find that: 
\begin{itemize}[itemsep=2pt,topsep=0pt,parsep=0pt]
    \item 
    User 1, belonging to the low-income student group, exhibits traits of a typical conservative and risk-averse investor in his pairwise item choices. 
    Both optimistic and neutral estimations can capture his characteristics {\color{black}well}. His utility value sharply increases under both estimations for relatively small returns ($\leq 10000$). 
    {\color{black} While his} marginal utility {\color{black}becomes} relatively small, even approaching zero, for {\color{black} relative} larger returns.  
    This indicates that the user gains significant satisfaction from wealth increase within lower returns and is less sensitive to changes in wealth for higher returns, aligning with his income level.  
    \item User 2 belongs to the middle-income middle-aged blue-collar group and has some savings.  
    Compared to User 1, User 2 exhibits characteristics of prudence and risk-averse, but with a decreasing risk-averse level. 
    \item User 3 belongs to the high-income young white-collar group and has some investment experience. 
    He essentially 
    {\color{black}has}
    a linear utility function, indicating {\color{black} an approximately risk-neutral attitude}, which is consistent with his choices on the pairwise items (always choosing the item with 
    {\color{black}the}
    higher expected return). 
    \item {\color{black}From RRA and ARA in Table \ref{tab:user-RRA}, we find that User 1 exhibits an overall higher degree of risk aversion than User 2, and User 2's degree of risk aversion is overall greater than that of User 3. In comparison, User 1 demonstrates higher risk aversion at low returns, while User 3 exhibits risk aversion only at high returns.}
\end{itemize}
{\color{black}All the observations are consistent to the personal characteristics of three participants, particularly the income level and investment budgets, shown in Table \ref{tab:user information}.
}

{\color{black}
\subsection{Application to human users: Analysis of risk aversion from elicited utility within population}\label{sec-real-port}
To quantify users' attitudes towards risk, akin to the concept of Gini coefficient in economics, we consider a measure, also named "Gini coefficient", based on elicited nominal utility functions to assess their degrees of risk aversion. 
Given a utility function $u(y)$, the Gini coefficient is defined as the area between functions $u(y)$ and $\frac{y}{\bar{y}_N}$ divided by $\frac{\bar{y}_N}{2}$:  
$${\rm Gini}(u(y)) 
=\frac{2\int_0^{\bar{y}_N} \left[u(y)-\frac{y}{\bar{y}_N}\right] \,dy}{\bar{y}_N}.
$$
The Gini coefficient of a user ranges between [0, 1]. A value closer to 1 indicates a higher degree of risk aversion, while a value closer to 0 suggests a tendency toward risk neutrality. 

In this subsection, we invite 60 participants with diverse characteristics to complete the same static preference questionnaire as that in subsection \ref{sec-real-utility}, which consists of the same 8 pairwise items from $\mathbb{I}^{10}$: $(I_1, I_7)$, $(I_4, I_2)$, $(I_2, I_5)$, $(I_6, I_{10})$, $(I_3, I_9)$, $(I_4, I_8)$, $(I_1, I_9)$ and $(I_3, I_{10})$. 
Then we elicit three nominal utility functions ($u^P$, $u^O$, $u^N$) and calculate three corresponding Gini coefficients for each participants. 
Given the personal information collected from these 60 participants additionally, we then conduct an analysis on how gender, age, and income levels influence individuals' levels of risk aversion, with relevant statistics shown in Table \ref{tab:gini}.

{\color{black}Moreover, we carry out the ANOVA (Analysis of Variance) test to make pairwise comparisons to the mean values of two segmented population under each characteristic. 
The corresponding $F$-value and the $p$-value for each pair of populations with three nominal utility functions are shown in Table \ref{tab:gini}. We set the significance level to $0.05$, thus a p-value less than $0.05$ indicates a significant difference in means within the pair.}

\begin{table}[htbp]
  \centering
  \caption{{\color{black}Statistics of Gini coefficients among 60 participants with different individual characteristics}}
\fontsize{8}{12}{
  \resizebox{\textwidth}{!}{\color{black}
    \begin{tabular}{c|c|c|cccc|cccc|cccc}
    \toprule
    \multicolumn{1}{r}{} &       &       & \multicolumn{4}{c|}{Gini($u^P$)} & \multicolumn{4}{c|}{Gini($u^O$)} & \multicolumn{4}{c}{Gini($u^N$)} \\
\cmidrule{3-15}    \multicolumn{1}{r}{} &       & Count & Mean  & Variance & $F$-value & $p$-value & Mean  & Variance & $F$-value & $p$-value & Mean  & Variance & $F$-value & $p$-value \\
    \midrule
    \multirow{2}[2]{*}{Gender} & Female & 30    & 0.0171  & 0.0004  & \multirow{2}[2]{*}{3.8488 } & \multirow{2}[2]{*}{0.0546 } & 0.9620  & 0.0020  & \multirow{2}[2]{*}{6.4351 } & \multirow{2}[2]{*}{0.0139 } & 0.4896  & 0.0004  & \multirow{2}[2]{*}{6.9348 } & \multirow{2}[2]{*}{0.0108 } \\
          & Male  & 30    & 0.0097  & 0.0000  &       &       & 0.8141  & 0.1000  &       &       & 0.4119  & 0.0257  &       &  \\
    \midrule
    \multirow{2}[2]{*}{Age} & Below 35 years & 34    & 0.0149  & 0.0003  & \multirow{2}[2]{*}{0.7427 } & \multirow{2}[2]{*}{0.3923 } & 0.9071  & 0.0414  & \multirow{2}[2]{*}{0.5063 } & \multirow{2}[2]{*}{0.4796 } & 0.4610  & 0.0107  & \multirow{2}[2]{*}{0.5700 } & \multirow{2}[2]{*}{0.4533 } \\
          & Above 35 years & 26    & 0.0115  & 0.0001  &       &       & 0.8631  & 0.0757  &       &       & 0.4373  & 0.0195  &       &  \\
    \midrule
    \multirow{2}[2]{*}{Monthly income} & Below 5000 & 27    & 0.0146  & 0.0002  & \multirow{2}[2]{*}{0.2841 } & \multirow{2}[2]{*}{0.5960 } & 0.9511  & 0.0020  & \multirow{2}[2]{*}{3.6605 } & \multirow{2}[2]{*}{0.0607 } & 0.4828  & 0.0005  & \multirow{2}[2]{*}{3.6803 } & \multirow{2}[2]{*}{0.0600 } \\
          & Above 5000 & 33    & 0.0125  & 0.0002  &       &       & 0.8365  & 0.0949  &       &       & 0.4245  & 0.0245  &       &  \\
    \bottomrule
    \end{tabular}}}
  \label{tab:gini}%
\end{table}%


The following conclusions can be derived from Table \ref{tab:gini}: 
\begin{itemize}
\item The mean values of Gini coefficients in the female group (aged under 35 or with income below 5000, respectively) are higher than those in the male group (aged above 35 or with income above 5000, respectively) under all three nominal utility functions. These results suggest overall trends where, on average, females, younger individuals, and lower-income people tend to be more risk-averse than males, older individuals, and higher-income people.
\item  The results of the ANOVA tests validate the observation. The differences in risk-aversion levels 
caused by gender are 
most significant and the differences caused by income are relatively significant. 
However, statistical tests do not support the 
the difference in age.
Excluding the age factor, the neutral or optimistic nominal utility functions provide more significant discrimination over the gender or income factor, compared to the pessimistic estimation.


\end{itemize}
}

\section{Conclusion}\label{sec-conclusion}
In this paper, we propose a robo-advisor system to provide personalized investment portfolio recommendations to users based on their individual risk preferences. 
The system first generates a static preference questionnaire consisting of several pairwise items to obtain pairwise comparison 
{\color{black}preferences}
of the user. 
Then the preference elicitation process estimates three nominal utility functions in a piecewise linear form, which {\color{black}are} consistent with the users' pairwise comparisons.  
Finally, the system recommends a portfolio to the user by maximizing her/his nominal utility. 

To verify the proposed robo-advisor system, 
we conduct some simulation tests and {\color{black} human experiments}. 
The results illustrate that 
the three nominal utility functions (pessimistic, optimistic, neutral) elicited 
converge to the preset true utility function as the number of questions increases, 
and the SPQ method exhibits a more efficient elicitation {\color{black} performance} compared to 
randomly generated questionnaires. 
Moreover, the nominal utility functions elicited for the three {\color{black} human} participants match their distinct personal characteristics well 
{\color{black}
and statistics observed from the risk aversion coefficients of the elicited nominal utility of the 60 human participants align with the corresponding group characteristics.
}

{\color{black}The static preference questionnaire generation method utilized in this paper is based on prior information, specifically the historical rating data. Exploring methods that select pairwise items without additional prior information such as \cite{szadoczki2022filling} is a promising topic.}
This proposed robo-advisor framework 
only involves a single interaction with the user through a static questionnaire, 
which 
may fall short in maintaining long-term contact with users.
Thus it's promising 
to develop a dynamic utility updating process that allows multiple interactions with users,  
{\color{black}
facilitating real-time adjustments of the utility function in response to evolving user preferences and requiring fewer questions to achieve the same level of reliability. For instance, the Swiss-system, which has been identified as an efficacious tournament design \cite{sziklai2022efficacy}, might be included.}
Moreover, {\color{black} human} users sometimes exhibit contradictory attitudes toward risk when answering questionnaires, 
which poses significant challenges in eliciting feasible nominal utility functions consistent with their preferences. 
Future 
efforts could focus on implementing additional mechanisms to identify conflicting responses. 
For instance, categorizing user answers based on different psychological states or employing incremental learning methods may be effective solutions. 

\footnotesize



\bibliography{ref}

\newpage
\appendix
\normalsize
\section{Test setting}
\begin{table}[htbp]
\setlength{\extrarowheight}{4pt}
  \centering
  \caption{Item set $\mathbb{I}^{20}$}\label{20-item set}
    \resizebox{\textwidth}{!}{   
    \begin{tabular}{l l l}
    \toprule
      Item & Description & Probability representation\\ 
      \midrule
    1 & \text{\textyen}100 in cash & $\mathbb{P}(I_1=100)=1$ \\
    2 & A 90\% chance at winning \text{\textyen}500 
    & $\mathbb{P}(I_2=500)=0.9,\ \mathbb{P}(I_2=0)=0.1$ \\
    3 & An 80\% chance at winning \text{\textyen}1,000 
    & $\mathbb{P}(I_3=1000)=0.8,\ \mathbb{P}(I_3=0)=0.2$ \\
    4 & A 60\% chance at winning \text{\textyen}3,000 
    & $\mathbb{P}(I_4=3000)=0.6,\ \mathbb{P}(I_4=0)=0.4$ \\
    5 & A 40\% chance at winning \text{\textyen}5,000 
    & $\mathbb{P}(I_5=5000)=0.4,\ \mathbb{P}(I_5=0)=0.6$ \\
    6 & A 20\% chance at winning \text{\textyen}10,000 
    & $\mathbb{P}(I_6=10000)=0.2,\ \mathbb{P}(I_6=0)=0.8$ \\
    7 & A 5\% chance at winning \text{\textyen}50,000 
    & $\mathbb{P}(I_7=50000)=0.05,\ \mathbb{P}(I_7=0)=0.95$ \\
    8 & A 3\% chance at winning \text{\textyen}100,000 
    & $\mathbb{P}(I_8=100000)=0.03,\ \mathbb{P}(I_8=0)=0.97$ \\
    9 & A 1.5\% chance at winning \text{\textyen}200,000 
    & $\mathbb{P}(I_9=200000)=0.015,\ \mathbb{P}(I_9=0)=0.985$ \\
    10 & A 1\% chance at winning \text{\textyen}300,000 
    & $\mathbb{P}(I_{10}=300000)=0.01,\ \mathbb{P}(I_{10}=0)=0.99$ \\
    11 & A 0.8\% chance at winning \text{\textyen}400,000 
    & $\mathbb{P}(I_{11}=400000)=0.008,\ \mathbb{P}(I_{11}=0)=0.992$ \\
    12 & A 0.6\% chance at winning \text{\textyen}500,000 
    & $\mathbb{P}(I_{12}=500000)=0.006,\ \mathbb{P}(I_{12}=0)=0.994$ \\
    13 & \parbox[t]{8cm}{An 88\% chance at winning \text{\textyen}500 \\ [-0.1ex] and a 0.5\% chance at winning \text{\textyen}500,000}
    & \parbox[t]{8.3cm}{$\mathbb{P}(I_{13}=0)=0.115,\ \mathbb{P}(I_{13}=500)=0.88,\\ [-0.1ex]  \mathbb{P}(I_{13}=500000)=0.005$}\\
    14 & \parbox[t]{8cm}{A 70\% chance at winning \text{\textyen}1,000 \\[-0.1ex]
    and a 0.6\% chance at winning \text{\textyen}400,000} 
    & \parbox[t]{8.3cm}{$\mathbb{P}(I_{14}=0)=0.294,\ \mathbb{P}(I_{14}=1000)=0.7,\\ [-0.1ex]  \mathbb{P}(I_{14}=400000)=0.006$} \\
    15 & \parbox[t]{8cm}{A 40\% chance at winning \text{\textyen}3,000 \\[-0.1ex]
    and a 0.6\% chance at winning \text{\textyen}300,000} 
    & \parbox[t]{8.3cm}{$\mathbb{P}(I_{15}=0)=0.594,\ \mathbb{P}(I_{15}=3000)=0.4,\\[-0.1ex]  \mathbb{P}(I_{15}=300000)=0.006$} \\
    16 & \parbox[t]{8cm}{A 25\% chance at winning \text{\textyen}7,000 \\[-0.1ex]
    and a 1\% chance at winning \text{\textyen}200,000} 
    & \parbox[t]{8.3cm}{$\mathbb{P}(I_{16}=0)=0.74,\ \mathbb{P}(I_{16}=7000)=0.25,\\[-0.1ex]  \mathbb{P}(I_{16}=200000)=0.01$} \\
    17 & \parbox[t]{8cm}{A 10\% chance at winning \text{\textyen}10,000 \\[-0.1ex]
    and a 2.5\% chance at winning \text{\textyen}100,000} 
    & \parbox[t]{8.3cm}{$\mathbb{P}(I_{17}=0)=0.875,\ \mathbb{P}(I_{17}=10000)=0.1,\\[-0.1ex]  \mathbb{P}(I_{17}=100000)=0.025$} \\
    18 & \parbox[t]{8cm}{A 20\% chance at winning \text{\textyen}7,000 \\[-0.1ex]
    and a 3\% chance at winning \text{\textyen}75,000} 
    & \parbox[t]{8.3cm}{$\mathbb{P}(I_{18}=0)=0.77,\ \mathbb{P}(I_{18}=7000)=0.2,\\[-0.1ex]  \mathbb{P}(I_{18}=75000)=0.03$} \\
    19 & \parbox[t]{8cm}{A 30\% chance at winning \text{\textyen}5,000 \\[-0.1ex]
    and a 3.5\% chance at winning \text{\textyen}50,000} 
    & \parbox[t]{8.3cm}{$\mathbb{P}(I_{19}=0)=0.665,\ \mathbb{P}(I_{19}=5000)=0.3,\\[-0.1ex]  \mathbb{P}(I_{19}=50000)=0.035$} \\
    20 & \parbox[t]{8cm}{A 90\% chance at winning \text{\textyen}100\\[-0.1ex]
    and a 10\% chance at winning \text{\textyen}25,000} 
    & $\mathbb{P}(I_{20}=100)=0.9,\ \mathbb{P}(I_{20}=25000)=0.1$\\
    \bottomrule
    \end{tabular}}
\end{table}

\end{document}